\documentclass[10pt]{amsart}
\allowdisplaybreaks[4]
\usepackage{amssymb,color}
\usepackage[square,compress,comma, numbers,sort]{natbib}

\usepackage[colorlinks=true, citecolor=blue, linkcolor=blue]{hyperref}
\usepackage{amsthm}
\usepackage{mathtools}

\newcommand{\E}[1]{\mathbb{E}\left\{ #1\right\}}

\newcommand{\pk}[1]{\mathbb{P} \left\{ #1 \right\} }

\def\Z{\mathbb{Z}}

\def\MB{\mathcal{B}}

\def\medt{\eta(d\vk{t})}
\def\medtl{\eta(d\vk{t})}

\def\sukt{\sigma_{\xi_{u,j}}(\vk{t})}

\def\suk0{\sigma_{\xi_{u,j}}(\vk{0})}

\def\ruk{\rho_{u,j}}
\def\vkt0{(\vk{t},\vk{0})}
\def\vks0{(\vk{s},\vk{0})}

\def\iukxz{\mathcal{I}_{u,j}(y;x)}
\def\KK{E}

\definecolor{c20}{rgb}{0.,0.7,0.}
\definecolor{c30}{rgb}{0.,0.,1.}
\definecolor{c40}{rgb}{1,0.1,0.7}
\definecolor{c50}{rgb}{1,0,0}
\definecolor{c60}{rgb}{0,0.9,0.1}

\newcommand{\vk}[1]{{#1}}

\newcommand{\abs}[1]{\left| #1 \right|}
\newcommand{\ABs}[1]{ \biggl \lvert #1 \biggr \rvert}
\newcommand{\norm}[1]{\lVert #1 \rVert}

\newcommand{\R}{\mathbb{R}}

\newcommand{\N}{\mathbb{N}}
\newcommand{\inr}{\in \R}

\newcommand{\ldot}{,\ldots,}

\newcommand{\limit}[1]{\lim_{#1 \to   \infty}}

\def\CC{\mathbb{Q}}

\newcommand{\BQN}{\begin{eqnarray}}
\newcommand{\EQN}{\end{eqnarray}}
\newcommand{\BQNY}{\begin{eqnarray*}}
\newcommand{\EQNY}{\end{eqnarray*}}

\def\cEE#1{\textcolor{black}{#1}}

\def\K1#1{\textcolor{cyan}{#1}}
\def\K1#1{#1}

\newcommand{\BS}{\begin{sat}}
\newcommand{\ES}{\end{sat}}
\newcommand{\BT}{\begin{theo}}
\newcommand{\ET}{\end{theo}}
\newcommand{\BK}{\begin{korr}}
\newcommand{\EK}{\end{korr}}

\newcommand{\BD}{\begin{de}}
\newcommand{\ED}{\end{de}}

\newcommand{\BIT}{\begin{itemize}}
\newcommand{\EIT}{\end{itemize}}
\newcommand{\BDI}{\begin{description}}
\newcommand{\EDI}{\end{description}}

\newcommand{\BRM}{\begin{remark}}
\newcommand{\ERM}{\end{remark}}

\newcommand{\BEL}{\begin{lem}}
\newcommand{\EEL}{\end{lem}}

\newtheorem{theo}{Theorem}[section]
\newtheorem{sat}[theo]{Proposition}
\newtheorem{de}[theo]{Definition}
\newtheorem{lem}[theo]{Lemma}

\newtheorem{korr}[theo]{Corollary}
\newtheorem{remark}[theo]{Remark}

\newcommand{\prooftheo}[1]{\textbf{Proof of Theorem} \ref{#1}}

\newcommand{\prooflem}[1]{\textbf{Proof of Lemma} \ref{#1}}

\newcommand{\COM}[1]{}

\newcommand{\QED}{\hfill $\Box$ \\}

\def\vp{\varepsilon}

\def\rw{\rightarrow}

\def\IF{\infty}

\def\CC{\mathbb{C}}
\def\LT{\left}
\def\RT{\right}

\def\KD#1{#1}

\def\EH#1{\textcolor{c40}{#1}}
\def\EH#1{#1}

\def\rd#1{\textcolor{red}{#1}}

\def\rd#1{#1}

\def\I{\mathbb{I}}
\def\rrd#1{\textcolor{black}{#1}}

\begin{document}

\title
{ Sojourn times of Gaussian related random fields}

\author{Krzysztof D\c{e}bicki}
\address{Krzysztof D\c{e}bicki, Mathematical Institute, University of Wroc\l aw, pl. Grunwaldzki 2/4, 50-384 Wroc\l aw, Poland}
\email{Krzysztof.Debicki@math.uni.wroc.pl}

\author{Enkelejd  Hashorva}
	\address{Enkelejd Hashorva, Department of Actuarial Science,
		University of Lausanne,
		UNIL-Dorigny, 1015 Lausanne, Switzerland
	}
	\email{Enkelejd.Hashorva@unil.ch}

\author{Peng Liu}
\address{Peng Liu, Department of Mathematical Sciences, University of Essex, Colchester, UK}
\email{peng.liu@essex.ac.uk}
	
\author{Zbigniew Michna}
\address{Zbigniew Michna, Department of Logistics, Wroc\l aw University of Economics and Business, Poland}
\email{Zbigniew.Michna@ue.wroc.pl}

\bigskip

\date{\today}
 \maketitle
\bigskip
{
	{\bf Abstract:} This paper is concerned with the asymptotic analysis of sojourn times of random fields with continuous
sample paths.
Under a very general framework we show that there is an interesting relationship
between tail asymptotics of sojourn times and that of supremum. Moreover, we establish the  uniform double-sum method
to derive the tail asymptotics of sojourn times.  In the literature, based on the pioneering research
of S. Berman the sojourn times have been utilised to derive the tail asymptotics of supremum of Gaussian processes.
In this paper we show that the opposite direction is even more fruitful, namely knowing the asymptotics of supremum o
f random processes and fields (in particular Gaussian) it is possible to establish the  asymptotics of their sojourn times.
We illustrate our findings considering i) two dimensional Gaussian random fields, ii) chi-process generated by stationary Gaussian processes and iii) stationary Gaussian queueing processes.


}

\section{Introduction \& First Result}\label{s.intro}

Let $X(t), t\in E$ be a random field  with
compact parameter set $E \subset \mathbb{R}^d,d\ge 1$ and \cEE{almost surely} continuous sample paths.
For a given level $u\inr$ define the excursion set of $X$ above the level $u$ by
$$ A_u(X)\coloneqq \{t\in E: X(t) >u \}.$$
The probability that $A_u$ is not empty
$$\pk{A_u(X) \not=\emptyset}= \pk{\exists t\in E: X(t)>u}= \pk{\sup_{t\in E}X(t)>u}=:p_u$$
is widely studied in the literature under the asymptotic regime $u\to \infty$, and the assumption that $X$ has marginals with infinite upper endpoint; see, e.g., \cite{Pit96, AdlerTaylor} for $X$ being Gaussian processes and related random fields.\\
Define the Lebesgue volume of $A_u(X)$ by

\def\vv{v }
\def\Veta{Vol}
$$ \Veta(A_u(X))=\int_{ \KK} \mathbb{I}(  X(t) >u)dt.$$
For specific cases, commonly $d=1$ and $X$ is stationary, asymptotic results as $u\to \IF$ are also known  for the probability that the volume of the excursion set (occupation time or sojourn time) exceeds $\vv(u)z$, i.e., approximations of
$$ r_u(z)\coloneqq\pk{\Veta(A_u(X)) > \vv(u)z}, \quad u\to \IF$$
for some specific positive scale function $\vv$ and $z\ge 0$ are available, see the seminal contribution \cite{Berman82}.

 The non-stationary case has been considered in \cite{Berman85, Berman87}. See also \cite{Berman92} for the comprehensive introduction of extremes of  sojourns for Gaussian processes.\\
In this contribution we are mainly interested in the formalisation of the uniform double-sum method for sojourns of random processes and fields focusing on the multidimensional case $d\ge 2$, for which no asymptotic results for  $r_u(z)$ are available in the literature.\\
 The first  question of our study is whether we can determine a  positive scaling functions $\vv(u),u>0$ and some survival function $\bar F$ such that
\begin{equation}\label{shti} \limit{u} \pk{ \Veta(A_u(X)) >\vv(u)z\Bigl\lvert\Veta(A_u(X))>0}=
	 	\limit{u} \pk{ Vol(A_u(X)) >\vv(u)z\Bigl\lvert \sup_{t\in E} X(t)>u}= \EH{\bar F}(z)
	\end{equation}
is valid for all $z\geq 0$. \\
If \eqref{shti} holds for some $z$ positive such that $\bar F(z) >0$ 
the asymptotics of $r_u(z)$ is proportional to that of
$p_u$, i.e.,
$$
 r_u(z) \sim \bar F(z) p_u, \quad u\to \IF.
$$
  Here $a(t)\sim b(t)$ means asymptotic equivalence of two real-valued functions $a(t)$ and $ b(t)$ when the argument $t$ tends to infinity or zero.
For a given index set $K$ we write $\sharp K$ for the cardinality of  $K$.

The following theorem states tractable conditions that imply \eqref{shti} for $X$ as above and \cEE{$E=E_u$. In order to avoid repetition, all Gaussian processes hereafter are assume to have almost surely continuous sample paths.}
\def\Esu{E(u,n)}
\BT\label{th1} 
\cEE{Let $E_u, u>0$  be compact set of $\R^d$ such that $\limit{u}\pk{ \sup_{t\in E_u} X(t)>u}=0 $. }
Suppose that there exist collections of Lebesgue measurable
disjoint compact sets $ I_{k}(u,n), k\in K_{u,n}$ with $K_{u,n}$ non-empty countable index sets    such that
$$
E(u,n)\coloneqq \bigcup_{k\in K_{u,n}}I_{k}(u,n)\subset \cEE{E_u},
$$
then \eqref{shti} holds with \cEE{$E=E_u$} if the following three conditions are satisfied:

A1) (Reduction to relevant sets)
$$ \lim_{n\rw\IF}\limsup_{u\rw\IF}\frac{\pk{ \sup_{t\in \cEE{E_u}\setminus\Esu} X(t)>u}}{\pk{ \sup_{t\in \cEE{\Esu}} X(t)>u}}=0.$$
A2) (Uniform single-sum approximation) There exists $v(u)>0$ and $\bar F_n, n\ge 1$ such that
\begin{equation}\label{Pickands}
\limit{u}\sup_{k\in K_{u,n}} \ABs{
		\frac{\pk{ Vol(\{ t\in I_{k}(u,n): X(t) >u \})> \vv(u)x}}
		 {\pk{ \sup_{ t\in I_{k}(u,n)} X(t)>u} }- \bar F_n(x)}= 0, \quad x\geq 0,\ {n\ge 1}
\end{equation}
and for all $x\geq 0$
\BQN\label{constant}
\bar F(x)\coloneqq \lim_{n\rw\IF} \bar F_n(x)\in (0,1].
\EQN
A3) (Double-sum negligibility) For all large $n$ and large $u$,  $\sharp K_{u,n}\geq 2$
and 
$$\lim_{n\rw\IF}\limsup_{u\rw\IF} \frac
{\sum_{ i \neq j, i,j\in K_{u,n}} \pk{ \sup_{t\in I_{i}(u,n)} X(t)>u,
\sup_{t\in I_j(u,n)} X(t)>u}}{\sum_{k\in K_{u,n}} \pk{ \sup_{t\in I_{k}(u,n)} X(t)>u}}= 0.
$$

\ET

For  $X(t),t\inr$ being a Gaussian process, \cite{KEP2016} shows that conditions A1)-A3) are satisfied under very general assumptions on $X$. From \cite{KEP2016},  we can formulate  some general conditions on $X$ that imply
\begin{equation} \limit{u}\sup_{k\in K_{u,n}}
		\ABs{\frac{ \pk{\sup_{ t\in I_{k}(u,n)} X(t)>u}}{ \Xi_{k}(u)}- C_{n}}= 0
		\label{C31}
	\end{equation}
for some known deterministic functions $\Xi_{k}(u)$, $k\in K_{u,n}$ and $C_n$  positive constants such that $\limit{n} C_n= C\in (0,\IF)$. In order to prove \eqref{Pickands} if \eqref{C31} holds, we shall prove that
	\begin{equation} \limit{u}\sup_{k\in K_{u,n}}
	\ABs{
		\frac{ \pk{ Vol(\{ t\in I_{k}(u,n): X(t) >u \})> \vv(u)x}}{ \Xi_{k}(u)}- D_{n}(x)}= 0, \quad
	\label{C32}
\end{equation}
where $D_n,\ n\ge1$ are deterministic functions such that $\limit{n} D_n(x)=D(x)>0$, $x\ge0$.
This then in turn implies that (\ref{constant}) holds with
$$
	 \bar F(x)= \frac{D(x)}{C}.
 $$
Note that in case that $D$ is continuous at $x=0$ we also expect that
$C=D(0)$
for all $z\ge 0$.\\

In the literature various results are known for supremum of functions of Gaussian vector processes,
for instance for chi-square processes, chaos of Gaussian processes, order statistics of Gaussian processes,
(see, e.g., \cite{Pitchi, Pit96, HJ2015,Bai19}) or reflected Gaussian processes
modelling a queueing process with Gaussian input (see, e.g.,
\cite{N1994,HP99,DE2002,Pit2001,Man03,HP2004,DI2005,Man07,KrzysPeng2015,KEP2018}).
In Section \ref{s.examples} we illustrate the applicability of Theorem \ref{th1}
by the analysis of three diverse families of stochastic processes:
1) Gaussian random fields (GRF's), 2) chi-processes and 3) reflected fractional Brownian motions.
For all this families of stochastic processes the available results
in the literature show that both A1) and A3)  hold
under quite general conditions; see Section 2.
Hence, in view of Theorem \ref{th1}, in order to get (\ref{shti})
it suffices to determine $\bar F$ in A2).

Except the above examples, our findings can also be applied to many other GRF's.
For instance, multi-dimensional GRF's with $d\geq 3$,
non-stationary chi-process or chi-square process, Gaussian chaos process, non-stationary Gaussian fluid queues and so on.
However, we shall not analyze these random processes or fields in this paper.

Brief organisation of the rest of the paper.
In Section \ref{s.Ber}
we introduce some notation and Berman-type constants
that
play the core role in the description of
$\bar F$.
In Section \ref{s.examples},
we provide examples that illustrate the derived in Theorem \ref{th1} technique for getting (\ref{shti}).
 Some technical lemmas are given in Section \ref{appen}; their proofs are deferred to Section
\ref{s.supl}. The proofs of the main contributions of this paper are presented in Section \ref{s.proofs}.

\section{Berman-type constants}\label{s.Ber}

We begin with the introduction  of the Berman-type constants for given independent fBm's
$B_{\alpha_i}(s),s\inr$ with Hurst index $\alpha_i/2\in (0,1]$, $i=1,2$.
 For given continuous functions $h_1,h_2$ set
 $$ W_{\alpha_1,\alpha_2,h_1,h_2}(t)\coloneqq \sum_{i=1}^2 (W_{\alpha_i}(t_i)  - h_i(t_i)), \quad t=(t_1,t_2)\inr^2,
 \quad W_{\alpha_i}(t_i)= \sqrt{2}B_{\alpha_i}(t_i) -\abs{t_i}^{\alpha_i}\,.
 $$
For simplicity, let $B_0(s)\equiv 0,s\inr$.
For 
$\alpha_i\in [0,2], i=1,2$,~ 
$x\geq 0$ and $E\subset \mathbb{R}^2$ a compact set, let
$$\mathcal{B}_{\alpha_1, \alpha_2}^{h_1, h_2}(x,E)=\int_{\mathbb{R}}\pk{\int_{ E}\I( W_{\alpha_1,\alpha_2,h_1,h_2}(t)>z) dt>x}e^zdz$$
and if the limit exists, define
$$\mathcal{B}_{\alpha_1, \alpha_2}^{b_1|t_1|^{\beta_1}, b_2|t_2|^{\beta_2}}(x)\coloneqq\lim_{S\rw\IF} \frac{\mathcal{B}_{\alpha_1, \alpha_2}^{b_1|t_1|^{\beta_1}, b_2|t_2|^{\beta_2}}(x, G(S,\alpha_1,\beta_1,\alpha_2,\beta_2)) }
{ S^{\mathbb{I}(\alpha_1<\beta_1)+\mathbb{I}(\alpha_2<\beta_2) }},$$
where
$$G(S,\alpha_1,\beta_1,\alpha_2,\beta_2)=\left\{\begin{array}{cc}
[0,S]^2,& \alpha_1<\beta_1, \alpha_2<\beta_2,\\
{[-S,S]\times [0,S]},&\alpha_1\geq \beta_1, \alpha_2<\beta_2,\\
 {[0,S]\times[-S,S]},& \alpha_1<\beta_1, \alpha_2\geq \beta_2,\\
{[-S,S]^2} & \alpha_1\geq \beta_1, \alpha_2\geq \beta_2.
\end{array}
\right.$$
We omit superscripts $h_i$'s if $h_1(s)=h_2(s)=0, s\in \mathbb{R}$ and then we put in our notation
 $\beta_1=\beta_2=\infty$ (this implies  that $\alpha_1<\beta_1$ and $\alpha_2<\beta_2$).
 Notice that for $x=0$, $\mathcal{B}_{\alpha_1, \alpha_2}^{h_1, h_2}(x)$ reduces to the classical Pickands or Piterbarg constants, see e.g., \cite{Pit96}. The one-dimensional Berman type constant is given by
$$\mathcal{B}_{\alpha}(x,[a,b])=\int_{\mathbb{R}}\pk{\int_{[a,b]}\I( W_{\alpha}(s)>z) ds>x}e^zdz$$
for $\alpha\in (0,2], a<b, a,b\in\mathbb{R}$, and
$$\mathcal{B}_{\alpha}(x)=\lim_{S\to\infty}\frac{\mathcal{B}_{\alpha}(x,[0,S])}{S}.$$
One can refer to \cite{KEZX17} and \cite{DLM18} for the existence and properties of
one-dimensional Berman constants.
For $x=0$, $\mathcal{H}_\alpha\coloneqq\mathcal{B}_{\alpha}(0)$ reduces to the classical Pickands constant; see, e.g., \cite{Pit96}.

The next lemma deals with properties of
\BQNY
\widehat{\mathcal{B}}_{\alpha_1,\dots,\alpha_m}\left(x, \prod_{i=1}^m[0,n_i]\right)\coloneqq\int_{\mathbb{R}}\pk{\int_{[0,n_1]}\I\left\{\sup_{t_i\in [0,n_i], i=2,\dots, m}\sum_{i=1}^{m}  W_{ \alpha_i}(t_i) >s\right\}dt_1>x}e^{s}ds
\EQNY
for $\alpha_i\in (0,2], i=1,\dots, m$ and $m\geq 1$.
\BEL\label{lemA0} For any $x\ge 0$, and $n_1>0$
\BQN
\label{Bermanconstant}\widehat{\mathcal{B}}_{\alpha_1,\dots,\alpha_m}(x,n_1)&\coloneqq&\lim_{n_i\rw\IF, i=2,\dots, m}\frac{\widehat{\mathcal{B}}_{\alpha_1,\dots,\alpha_m}\left(x, \prod_{i=1}^m[0,n_i]\right)}{\prod_{i=2}^{m}n_i} \notag \\
& =&
\prod_{i=2}^m \mathcal{H}_{\alpha_i}  \int_{\mathbb{R}}\pk{  \int_{ [0,n_1]} \I\left\{ W_{ \alpha_1}(t)>s \right\} dt>x} e^{s}ds \in (0,\infty)
\EQN
and
\BQN
\widehat{\mathcal{B}}_{\alpha_1,\dots,\alpha_m}(x)\coloneqq\lim_{n\rw\IF}\frac{\widehat{\mathcal{B}}_{\alpha_1,\dots,\alpha_m}(x,n)}{n}
= \mathcal{B}_{\alpha_1}(x)\prod_{i=2}^m \mathcal{H}_{\alpha_i}\in (0,\infty).
\EQN
\EEL
%
\begin{remark}
\label{lemma} The limits in (\ref{Bermanconstant}) are finite and positive
and $\widehat{\mathcal{B}}_{\alpha_1,\dots,\alpha_m}(x,n_1)$ is
a continuous function of $x$  over $[0,n_1)$ which follows from
the combination of Lemma \ref{lemA0} and Lemma 4.1 in \cite{DLM18}.
The claim of Lemma \ref{lemA0}  still holds if we replace $B_{\alpha_i}$ by $X_i$
being independent centered Gaussian processes with stationary increments
and variance function satisfying some regular conditions as e.g. in \cite{DE2002}.
\end{remark}

\section{Illustrating examples}\label{s.examples}
In this section
we shall apply  Theorem \ref{th1} to three classes of processes:
i) GRF's , ii) chi-process generated by {a stationary} Gaussian process  and
iii) stationary reflected fractional Brownian motions with drift.
\subsection{Sojourns of GRF's }

Although numerous results for the tail asymptotics of supremum of GRF's are available
for both stationary and non-stationary cases (see e.g., \cite{Pit96,Pit20}), sojourns
have not been treated so far in the literature. It follows from the available results
in the literature, that A1) holds under quite general conditions, for instance when the variance
function has a unique point of maximum and $X$ satisfies a global H\"older continuity condition, see e.g., \cite{Pit96}.
The main tool for proving A1) is the so-called Piterbarg inequality, see \cite{Pit96}[Thm 8.1] and the recent contribution
\cite{KEP2015}. Under some further weak assumptions on the variance/covariance function of $X$, also A3) has been
shown to hold for a wide collection of cases of interest, see \cite{Pit96,nonhomoANN}.
Thus, in light of Theorem \ref{th1}, in order to prove (\ref{shti}) for GRF's the main task is the
explicit calculation of $\bar F$.

\subsubsection{GRF's with constant variance}
First  we consider  $X$ being a centred  GRF with $Var(X(t))=1, t\in E\subset \mathbb{R}^2$ and the correlation function $r(t, s)$, $t,s\in \mathbb{R}^2$ satisfying
\BQN\label{cor}
1-r(t_1,t_2, s_1,s_2)\sim a_1|t_1-s_1|^{\alpha_1}+a_2|t_2-s_2|^{\alpha_2}, \quad (t_1,t_2), (s_1,s_2)\in E,  |t_i-s_i|\rw 0, i=1,2,
\EQN
with $a_i>0$ and $\alpha_i\in (0,2]$, $i=1,2$. Moreover,
\BQN\label{cor1}
r(t_1,t_2, s_1,s_2)<1, \quad (t_1,t_2), (s_1,s_2)\in E, (t_1,t_2)\neq (s_1,s_2).
\EQN
\cEE{
For notational simplicity  we shall consider $E=[0,T_1]\times [0,T_2]$, the results for general hypercubes in $\R^d$ follows with similar calculations. The case that $T_i=T_{i,u},i=1,2$ depend on $u$ needs some extra care. $T_{i,u}$'s should not be too small, i.e.,
$$\limit{u} T_{i,u}u^{2/\alpha_i}=\IF, i=1,2.$$
 On the other side $T_{i,u}$'s cannot be too large too. If the GRF  is stationary, then for some $\beta\in (0,1)$ we should require that $$\limit{u}T_{1,u}T_{2,u}e^{- \beta u^2/2}=0.$$
  In the more complex situation that we are looking at below the existence of $\beta$ is not clear.
We suppress the discussion for long intervals in order to avoid further complications. }
\BS\label{THSTA} Let $X(t), t\in E=[0,T_1]\times [0,T_2]$ be a centred  GRF
which satisfies (\ref{cor}) and (\ref{cor1}) and assume that
$\vv(u)=a_1^{-1/\alpha_1}a_2^{-1/\alpha_2} u^{-2/\alpha_1-2/\alpha_2}$. Then for all $x\ge0$
\BQNY 
\lim_{u\rw\IF}\pk{\int_{E}\I(X(t)>u)dt>\vv(u)x \Bigl\lvert \sup_{t\in E} X(t)>u}
= 
\frac{\mathcal{B}_{\alpha_1, \alpha_2}(x)}{\mathcal{B}_{\alpha_1, \alpha_2}(0)}. 
\EQNY
\ES

\subsubsection{GRF's with non-constant variance}
Denote by $\sigma(t)=\sqrt{Var(X(t))}$ and
assume that \\
$t^*={(t^*_1, t^*_2)}\in E=[-T_1,T_1]\times [-T_2,T_2]$ is the inner point of $E$,
which is the unique point such that $\sigma(\EH{t^*})=\sup_{ t\in E}\sigma(t)=1$ satisfying
\BQN\label{var}
1-\sigma(t)\sim b_1|t_1-t_1^*|^{\beta_1}+b_2|t_2-t_2^*|^{\beta_2}, \quad t=(t_1,t_2)\in E,  \,{\norm{t-t^*} \rw 0},
\EQN
with $b_i>0, \beta_i>0$, $i=1,2$. Here $\norm{\cdot}$ denotes the Euclidean norm. Moreover, let
\BQN\label{cor2}
1-r(t, s)\sim a_1|t_1-s_1|^{\alpha_1}+a_2|t_2-s_2|^{\alpha_2}
\EQN
as $t,s\in E, \norm{t-t^*},  \norm{s-t^*} \rw 0$ with $a_i>0$ and $\alpha_i\in (0,2]$, $i=1,2$, $s=(s_1,s_2)$,
where $r(t,s)$ is the correlation function of the random field $X$.
In the notation below we interpret $\infty\cdot 0$ as $0$.
\BS\label{Onepoint}  If $X(t), t\in E$ is a centered GRF  which satisfies (\ref{var}) and (\ref{cor2}) and   $\vv(u)= \prod_{i=1}^2\left(a_i^{-1/\alpha_i^*} u^{-2/\min(\alpha_i,\beta_i)}\right)$
with $\alpha_i^*=\alpha_i \I(\alpha_i\leq \beta_i)+\IF \I(\alpha_i> \beta_i)$, then for all $x\ge0$ 
$$
\lim_{u\rw\IF}\pk{\int_{E}\I(X(t)>u)dt>\vv(u)x \Bigl\lvert \sup_{t\in E} X(t)>u}=
\frac{\mathcal{B}_{\hat{\alpha}_1, \hat{\alpha}_2}^{\bar a_1 b_1 \abs{t_1}^{\beta_1},
\bar a_2 b_2 \abs{t_2}^{\beta_2}} (x)}{\mathcal{B}_{\hat{\alpha}_1, \hat{\alpha}_2}^{\bar a_1 b_1 \abs{t_1}^{\beta_1},
\bar a_2 b_2 \abs{t_2}^{\beta_2}} (0)},$$
where
  $$\bar a_i=\left\{\begin{array}{cc}
  0& \alpha_i<\beta_i\\
  \frac{1}{a_i}&\alpha_i=\beta_i\\
  1& \alpha_i>\beta_i
  \end{array}\right. ,\quad \hat{\alpha}_i=\left\{\begin{array}{cc}
  \alpha_i& \alpha_i\leq \beta_i\\
  0&\alpha_i>\beta_i
  \end{array}\right., i=1,2.$$
 \ES

\subsection{Sojourns of chi-processes}
Let $X(t), t\in [0,T]$ be a centered stationary 
Gaussian process with unit variance and correlation function satisfying
\BQNY
1-r(s,t)\sim a |t-s|^{\alpha}, \quad |s-t|\rw 0,
\EQNY
where $\alpha\in(0,2]$ and for all $s\neq t$, $s,t\in [0,T]$
\BQNY
r(s,t)<1.
\EQNY


Define the chi-process of degree $m\geq 1$ by
\BQN\label{chi}
\chi(t)\coloneqq\sqrt{\sum_{i=1}^m X_i^2(t)}, \quad t\in\mathbb{R},
\EQN
where $X_i, 1\leq i\leq m$ are iid copies of $X$. The exact asymptotics of
$\pk{\sup_{t\in [0,T]}\chi(t)>u}$ \cEE{has} been investigated  in \cite{Pitchi, Pit96, HJ2015}.
In the following theorem we 
consider 
the sojourn time of $\chi$.
\BS\label{chith} Let $\chi$ be defined as in (\ref{chi}). If $\vv(u)=a^{-1/\alpha} u^{-2/\alpha}$, then for all $x\geq 0$
$$\lim_{u\rw\IF}\pk{\int_{[0,T]}\I(\chi(t)>u)dt>\vv(u)x \Bigl\lvert \sup_{t\in [0,T]}\chi(t)>u}=
\frac{{\mathcal{B}}_{\alpha}(x)}{{\mathcal{B}}_{\alpha}(0)}.
$$
\ES
\subsection{Sojourns of stationary reflected fractional Brownian motion with drift}
Consider a stationary
reflected fractional Brownian motion with drift $Q(t),t\geq 0$, i.e.,
\BQNY
Q(t)\coloneqq\sup_{s\geq t} \left(B_{\alpha}(s)-B_{\alpha}(t)-c(s-t)\right),
\EQNY
where $B_\alpha$ is an fBm with Hurst parameter $\alpha/2\in (0,1)$ and $c\in (0,\IF)$.
Motivated by some applications of $Q(t)$ to queueing models,
the seminal paper \cite{HP99} studied the tail asymptotics of $Q(0)$.
Later on, \cite{Pit2001} considered the tail asymptotics of the supremum of $Q(t)$ over a time horizon.
Recently, the findings of Piterbarg have been extended to Gaussian processes with stationary increments
 \cite{KrzysPeng2015}. We consider next the case of fBm and note that a more general case of Gaussian processes with stationary increments can be also dealt with using results from \cite{KrzysPeng2015}.
In the following we consider \cEE{$E_u=[0,T_u]$}, where $T_u$ is a non-negative function of $u>0$.

\BS\label{queue} Let $\vv(u)=u^{\frac{2(\alpha-1)}{\alpha}}\left(\frac{\sqrt{2}(\tau^*)^\alpha}{1+c\tau^*}\right)^{2/\alpha}$
with $\tau^*=\frac{\alpha}{c(2-\alpha)}$ and $\alpha\in(0,2)$. \\
i) If $\limit{u} \frac{T_u}{v(u)} = T\in (0,\IF)$, then for $T>x\geq 0$
$$\lim_{u\rw\IF}\pk{\int_{[0,T_u]}\I(Q(t)>u)dt>\vv(u)x\Bigl\lvert \sup_{t\in [0,T_u]}Q(t)>u}=\frac{\mathcal{B}_{\alpha}(x, [0,T])}{\mathcal{B}_{\alpha}(0, [0,T])}. 
$$
ii) If $ \limit{u} \frac{T_u}{v(u)} =\IF$ and $T_u<e^{\beta u^{2-\alpha}}$ with $\beta\in \left(0, \left(\frac{1+c\tau^*}{\sqrt{2}(\tau^*)^{\alpha/2}}\right)^2\right)$, then for all $x\geq 0$
$$\lim_{u\rw\IF}\pk{\int_{[0,T_u]}\I(Q(t)>u)dt>\vv(u)x \Bigl\lvert \sup_{t\in [0,T_u]}Q(t)>u}=
\frac{{\mathcal{B}}_{\alpha}(x)}{{\mathcal{B}}_{\alpha}(0)}.
$$
\ES
\BRM 1) Note that $\lim_{u\rw\IF}v(u)=\IF$ for $\alpha>1$,
and $\lim_{u\rw\IF}v(u)=0$ for $\alpha<1$.

2) Conclusion in i)  of Proposition \ref{queue} still holds for $x>T$ since both sides in the equality of i) are $0$. However, it becomes tricky for the case $T=x$. We consider two special cases for $T=x$. If $T=x$ and $T_u\leq xv(u)$ for $u$ sufficiently large, then both sides in the equality of i) are  $0$. If $T=x$ and $T_u>xv(u)$ for sufficiently large $u$, we get, as $u\to\infty$
$$\pk{\int_{[0,T_u]}\I(Q(t)>u)dt>\vv(u)x\Bigl\lvert \sup_{t\in [0,T_u]}Q(t)>u}\sim \frac{\pk{\inf_{t\in [0,T_u]}Q(t)>u}}{\pk{\sup_{t\in [0,T_u]}Q(t)>u}}.$$
Combining the above two cases for $T=x$, we conclude that the limit for $T=x$ generally does not exist.
\ERM


\section{Auxiliary lemmas}\label{appen}
In this section we collect some lemmas that play important, although mostly technical role in the proofs of
results given in Sections \ref{s.intro}-\ref{s.examples}.
Their proofs are deferred to Section \ref{s.supl}. We begin with a
lemma {which is} an extension of Theorem 2.1 from \cite{KEP2016}.
Suppose that
for a compact $d-$dimensional hyperrectangle $K\subset\R^d$ we have
$$I_{k}(u,n)=\{t_{u,n,k}+(v_1(u)t_{{1}}, \dots v_d(u)t_d): \, t\in K\},$$
where $v_i(u)>0$, $i=1,\ldots,d$ and $t=(t_1,\ldots,t_d)\inr^d$.
Then, by  transforming time, we have
\BQNY\mathbb{P}\left(Vol(\{ t\in I_{k}(u,n): X(t) >u \})>v(u)z \right)&=& \mathbb{P}\left(\int_{ I_{k}(u,n)}\mathbb{I}(X(t)>u)dt>v(u)z\right)\\
&=&\mathbb{P}\left(\int_{\KD{K}} \mathbb{I} (X(t_{u,n,k}+(v_1(u)t_{{1}}, \dots v_d(u)t_d) )> u) dt>z\right),
\EQNY
where  $v(u)=\prod_{i=1}^dv_i(u)$.\\
\cEE{Motivated by these calculations, we  consider next}
$\xi_{u,j}(t),  t\in \KK_1, \ j\in S_u, \ {u\geq 0}$
a family of centered GRF's with continuous sample paths and variance function $\sigma_{u,j}^2$. \\
\cEE{Suppose in the following that $S_u$ is a countable set for all $u$ large.}\\
For simplicity in the following we assume that $\vk{0}\in E_1$. For a random variable $Z$, we set $\overline{Z}=\frac{Z}{\sqrt{Var(Z)}}$ if $Var(Z)>0$.

 We introduce next three assumptions:
 \begin{itemize}
 	\item[{\bf C0}:]   $\{g_{u,j}, j\in S_u\}$ is a sequence of deterministic functions of $u$ satisfying
 	\BQNY
 	\lim_{u\to\IF}\inf_{j\in S_u}g_{u,j}=\IF.
 	\EQNY
 	\item[{\bf C1}:] $Var(\xi_{u,j}(\vk{0}))=1$ for all large $u$ and any $j\in S_u$ and
 	there exists some bounded continuous function $h$ on $ \KK_1$ such that
 	\BQNY\label{assump-cova-field}
 	\lim_{u\to\IF}\sup_{\vk{s}\in E_1 ,j\in S_u}\left|g_{u,j}^2\LT( 1-\sigma_{u,j}(\vk{s}) \RT) - h(\vk{s})\right| =0.
 	\EQNY	
 	\item[{\bf C2}:] There exists a centered GRF  $\zeta(\vk{s}),\vk{s}\in \mathbb{R}^{d}$ with a.s. continuous sample paths such that 
\BQN\label{C21}
 	\lim_{u\to\IF}\sup_{s, s'\in E_1, j\in S_u}\abs{g_{u,j}^2\big(Var(\overline\xi_{u,j}(\vk{s})-\overline\xi_{u,j}(\vk{s}'))\big) - 2Var(\zeta(s)-\zeta(s'))} =0.
\EQN
 	\item[{\bf C3}:] There exist positive constants $C, \nu, u_0$ such that
 	\BQNY\label{assump-holder-field}
 	\sup_{j\in S_u} g_{u,j}^2Var(\overline\xi_{u,j}(\vk{s})-\overline\xi_{u,j}(\vk{s}')) \leq C \norm{\vk{s}-\vk{s}'}^\nu
 	\EQNY
 	holds for all $\vk{s},\vk{s}'\in \KK_1 , u\geq u_0$.
 \end{itemize}
Denote by $C(E_{i}), i=1,2$ \cEE{the Banach space of  all continuous functions $f: E_i \mapsto \R$},  with $ E_{i}\subset\mathbb{R}^{d_i}, d_i\geq 1, i=1,2$ being  compact rectangles \cEE{equipped with the sup-norm.}

Let $\Gamma: C(E_{1})\rw C(E_{2})$ be a continuous functional satisfying
\\
{\bf F1}: For any $f\in C(E_{1})$, and $a>0, b\in\mathbb{R}$, $\Gamma(af+b)=a\Gamma(f)+b$;\\
{\bf F2}: There exists $c>0$ such that
 $$\sup_{t\in E_{2}}\Gamma(f)(t)\leq c\sup_{s\in E_{1}}f(s), \ \ \forall f\in C(E_{1}).$$
Hereafter, $Q_i, i\in \N$ are some positive constants which might be different from line to line and $f(u,n)\sim g(u), u\to\infty, n\to\infty$ means that $$\lim_{n\to\infty}\lim_{u\to\infty}\frac{f(u,n)}{g(u)}=1.$$
\BEL\label{the-weak-conv}
Let $\{\xi_{u,j}(\vk{s}),\vk{s}\in E_{1} ,j\in S_u, {u\geq 0}\}$ be a family of centered GRF's defined as above satisfying {\bf C0-C3} and let
$\Gamma$  satisfy  {\bf F1-F2}.
Let $\eta$ be a positive $\sigma$-finite measure on
	$E_2$ being equivalent with the Lebesgues measure on $E_2$. \cEE{If for  all} large
	$u$ and all $j \in S_u$
$$\pk{ \sup_{t\in E_{2}} \Gamma(\xi_{u,j})(t)> g_{u,j}} > 0,$$
then for all $x\in[0,\eta(E_{2}))$
\BQN\label{con-uni-con}
\lim_{u\to\IF}\sup_{j\in S_u} \ABs{ \frac{\pk{ \int_{ E_{2} } \mathbb{I}\LT(\Gamma(\xi_{u,j})(\vk{t})>g_{u,j}\RT) \eta(dt) >x} } {\Psi(g_{u,j})} -\MB^{\Gamma, h,\eta}_{\zeta}( x, E_2) } = 0,
\EQN where $\Psi$ is the tail of the standard normal distribution and
\BQNY
\MB^{\Gamma, h,\eta}_{\zeta}( x, E_2)\coloneqq \int_\R  \pk{ \int_{ E_{2} } \mathbb{I}\big(\Gamma(\sqrt{2}\zeta-Var(\zeta)-h)(t) + y >0\big) \eta(dt) >x } e^{-y} dy
\EQNY
and  the constant $\MB^{\Gamma, h,\eta}_{\zeta}( x, E_2)$ is continuous \cEE{at  $x\in {(}0,\eta(E_{2}))$}.
\EEL


\BEL\label{l.1}
Let $x\geq 0$. Then
\\
\vspace{2mm}
(i) $\mathcal{B}_{\alpha_1,\alpha_2}(x)=\lim_{n\rw\IF}\frac{\mathcal{B}_{\alpha_1,\alpha_2}(x, [0,n]^2)}{n^2}\in (0,\IF),$
\\
\vspace{2mm}
%
%
(ii) $\lim_{n\rw\IF}\frac{\mathcal{B}_{\alpha_1,\alpha_2}^{a_1^{-1}b_1|t_1|^{\alpha_1},0}(x, [-n,n]\times[0,n])}{n}\in (0,\IF),
$
\\
\vspace{2mm}
%
(iii)
$
 \lim_{n\rw\IF} \mathcal{B}_{\alpha_1,\alpha_2}^{a_1^{-1}b_1|t_1|^{\alpha_1},a_2^{-1}b_2|t_2|^{\alpha_2}}(x, [-n,n]^2)\in (0,\IF).
$
%
\EEL


\section{Proofs}\label{s.proofs}

\subsection{\prooftheo{th1}}
\cEE{Let next $ A_u(X)\coloneqq \{t\in E_u: X(t) >u \}.$} For all $x\ge 0$ and all $u$ positive, since $v(u)$ is non-negative we have
\BQNY
\pi(u)&\coloneqq& \pk{ Vol(A_u(X)) >\vv(u)x\Bigl\lvert Vol(A_u(X))> 0}\\
&=& \pk{ Vol(A_u(X)) >\vv(u)x\Bigl\lvert \sup_{t\in E_u} X(t)>u}\\
&=&\frac{\pk{\int_{E}\mathbb{I}(X(t)>u)dt >\vv(u)x}}{\pk{\sup_{t\in \cEE{E_u}} X(t)>u}}
\EQNY
and further for all $n\ge 1$
\BQNY
\pi(u)&\geq& \frac{\pk{\int_{E(u,n)}\mathbb{I}(X(t)>u)dt >\vv(u)x}}{\pk{\sup_{t\in E(u,n)} X(t)>u}+\pk{\sup_{t\in \cEE{E_u}\setminus E(u,n)} X(t)>u}},\\
\pi(u)
&\leq& \frac{\pk{\int_{E(u,n)}\mathbb{I}(X(t)>u)dt >\vv(u)x}}{\pk{\sup_{t\in E(u,n)} X(t)>u}}+\frac{\pk{\sup_{t\in  \cEE{E_u} \setminus E(u,n)} X(t)>u}}{\pk{\sup_{t\in  \cEE{E(u,n)}} X(t)>u}}.
\EQNY
Applying {A1}, it follows that
\BQNY
\pi(u)\sim \frac{\pk{\int_{E(u,n)}\mathbb{I}(X(t)>u)dt >\vv(u)x}}{\pk{\sup_{t\in E(u,n)} X(t)>u}}=:\pi(u,n), \quad u\rw\IF, n\rw\IF.
\EQNY
For the case that $\sharp K_{u,n}=1$ for $u$ and $n$ sufficiently large,
the claim can be established straightforwardly by A2. Thus let us suppose that $\sharp K_{u,n}\geq 2$ for  $n$ and $u$ sufficiently large.
\cEE{In order to proceed we shall apply the standard scheme utilising Bonferroni inequality. Set therefore }
$$\Sigma_{u,n}\coloneqq\sum_{k\in K_{u,n}} \pk{ \sup_{t\in I_{k}(u,n)} X(t)>u},\quad \Sigma\Sigma_{u,n}\coloneqq\sum_{ i \neq j, i,j\in K_{u,n}} \pk{ \sup_{t\in I_i(u,n)} X(t)>u,
\sup_{t\in I_j(u,n))} X(t)>u}.$$
By the Bonferroni  inequality
\BQNY
\Sigma_{u,n}-\Sigma\Sigma_{u,n}\leq \pk{\sup_{t\in E(u,n)}X(t)>u}\leq \Sigma_{u,n}.
\EQNY
\cEE{The asymptotic behaviour of the probability of interest in the above inequality can be derived if the following two-step procedure is successful (which will work in our settings here). First we determine the exact asymptotics of the upper bound and then in a second step we  show that the correction in the lower bound is asymptotically negligible.\\
Now we want to apply the same idea for the sojourn functional, here the analysis is however more involved.  Observe first that } for any $u>0$
\BQNY
&&\pk{\int_{E(u,n)}\mathbb{I}(X(t)>u) dt >\vv(u)x}\\
&&\quad\leq\pk{\sum_{k\in K_{u,n}}\int_{I_{k}(u,n)}\mathbb{I}(X(t)>u) dt >\vv(u)x}\\
&&\quad \leq \pk{ \exists  k\in K_{u,n}, \int_{I_{k}(u,n)}\mathbb{I}(X(t)>u) dt >\vv(u)x}\\
&&\quad \quad + \pk{ \exists i, j\in K_{u,n}, i\neq j, \int_{I_i(u,n)}\mathbb{I}(X(t)>u) dt>0, \int_{I_j(u,n)}\mathbb{I}(X(t)>u) dt>0 }\\
&&\quad \leq \hat{\pi}(u,n)+\Sigma\Sigma_{u,n},
\EQNY
where
$$\hat{\pi}(u,n)=\sum_{k\in K_{u,n}}\pk{\int_{I_{k}(u,n)}\mathbb{I}(X(t)>u) dt >\vv(u)x}.$$
Using Bonferroni inequality again we have
\BQNY
\pk{\int_{E(u,n)}\mathbb{I}(X(t)>u) dt >\vv(u)x}&\geq& \pk{ \exists k\in K_{u,n}, \int_{I_{k}(u,n)}\mathbb{I}(X(t)>u) dt >\vv(u)x}\\
&\geq &\hat{\pi}(u,n)-\Sigma\Sigma_{u,n}.
\EQNY
 \cEE{The sojourn integral can then be approximated by $\hat{\pi}(u,n)$ if we show the correction in the lower bound is negligible. }
We have
 \BQNY
 \limsup_{u\rw\IF}\pi(u,n)&\leq&\limsup_{u\rw\IF}\frac{\hat{\pi}(u,n)+\Sigma\Sigma_{u,n}}{\Sigma_{u,n}-\Sigma\Sigma_{u,n}}= \limsup_{u\rw\IF}\frac{\hat{\pi}(u,n)}{\Sigma_{u,n}}\times \frac{1+\limsup_{u\rw\IF}\frac{\Sigma\Sigma_{u,n}}{\hat{\pi}(u,n)}}{1-\limsup_{u\rw\IF}\frac{\Sigma\Sigma_{u,n}}{\Sigma_{u,n}}},\\
 \liminf_{u\rw\IF}\pi(u,n)&\geq&\liminf_{u\rw\IF}\frac{\hat{\pi}(u,n)-\Sigma\Sigma_{u,n}}{\Sigma_{u,n}}=\liminf_{u\rw\IF}\frac{\hat{\pi}(u,n)}{\Sigma_{u,n}}-\limsup_{u\rw\IF}\frac{\Sigma\Sigma_{u,n}}{\Sigma_{u,n}}.
\EQNY
By (\ref{Pickands}) in A2 for any $n\geq 1$ and $x\geq 0$
\BQNY
\limsup_{u\rw\IF}\frac{\hat{\pi}(u,n)}{\Sigma_{u,n}}=\liminf_{u\rw\IF}\frac{\hat{\pi}(u,n)}{\Sigma_{u,n}}=\bar F_n(x)
\EQNY
implying
\BQN\label{barF}
\bar F_n(x)-\limsup_{u\rw\IF}\frac{\Sigma\Sigma_{u,n}}{\Sigma_{u,n}}\leq \liminf_{u\rw\IF}\pi(u,n)\leq \limsup_{u\rw\IF}\pi(u,n)\leq \bar  F_n(x)\times \frac{1+\limsup_{u\rw\IF}\frac{\Sigma\Sigma_{u,n}}{ \bar F_n(x)\Sigma_{u,n}}}{1-\limsup_{u\rw\IF}\frac{\Sigma\Sigma_{u,n}}{\Sigma_{u,n}}}.
\EQN
In view of A3, letting $n\rw\IF$ in the above inequalities we have that for $x\geq 0$
$$\lim_{n\rw\IF}\lim_{u\rw\IF}\pi(u,n)= \bar F(x)\in {(0,1]}.$$
This completes the proof. \QED

\subsection{\prooflem{lemA0}}
By the independence of $W_{\alpha_i}$'s for any positive $n_1 \ldot n_m$
\BQNY
\widehat{\mathcal{B}}_{\alpha_1,\dots,\alpha_m}\left(x, \prod_{i=1}^m[0,n_i]\right)&=&
\E{  \int_{\mathbb{R}}\I( \int_{[0,n_1]}\I\left\{\sup_{t_i\in [0,n_i], i=2,\dots, m}\sum_{i=1}^{m}  W_{ \alpha_i}(t_i)>s \right\}dt_1>x) e^{s}ds}\\
&=&\E{ e^{\sum_{i=2}^{m}  \sup_{t_i \in [0,n_i]}W_{ \alpha_i}(t_i)   }
	\int_{\mathbb{R}}\I( \int_{[0,n_1]}\I\left\{ W_{ \alpha_1}(t_1)>s \right\} dt_1 >x) e^{s} ds}\\
&=& \prod_{i=2}^m\E{ \sup_{t_i \in [0,n_i]} e^{  W_{ \alpha_i}(t_i)}   }
\int_{\mathbb{R}}\pk{  \int_{ [0,n_1]} \I\left\{ W_{ \alpha_1}(t_1)>s \right\} dt_1>x} e^{s}ds.
\EQNY
Hence the  claim follows by the definition of Pickands and Berman constants.\QED

\subsection{Proof of Proposition \ref{THSTA}}
The proof will be established by
checking that  A1-A3 in Theorem \ref{th1} are satisfied.
We begin with the introduction of partition
$$
I_{k_1,k_2}(u,n)=\prod_{i=1}^2[a_i^{-1/\alpha_i}u^{-2/\alpha_i}k_in, a_i^{-1/\alpha_i}u^{-2/\alpha_i}(k_i+1)n],
$$
for $$ 0\leq k_i\leq [T_ia_i^{1/\alpha_i}u^{2/\alpha_i}n^{-1}]-1=:N_{i}(u,n),\ i=1,2.$$
Let
$$K_{u,n}=\{(k_1,k_2): 0\leq k_1\leq N_1(u,n), 0\leq k_2\leq N_2(u,n)\}$$
and
$E(u,n)=\bigcup_{(k_1,k_2)\in K_{u,n}} I_{k_1,k_2}(u,n)$.
Then $E(u,n)\subset E$.

{\it \underline{Condition A1.}}
It follows straightforwardly
from Lemma 7.1 in \cite{Pit96} that
\BQN\label{Stationary}\pk{ \sup_{t\in E} X(t)>u}
\sim \sum_{0\leq k_i\leq N_i(u,n), i=1,2}\pk{ \sup_{t\in I_{k_1,k_2}(u,{n})} X(t)>u}
, \quad u\rw\IF,
n\to\infty,
\EQN
which implies that
condition A1 holds.\\

{\it \underline{Condition A2.}}  Let for $t=(t_1, t_2)$
$$\xi_{u,n,k_1,k_2}(t)=X(a_1^{-1/\alpha_1} u^{-2/\alpha_1}(k_1n+t_1), a_2^{-1/\alpha_2} u^{-2/\alpha_2}(k_2 n+t_2)), \quad \vv(u)=a_1^{-1/\alpha_1}a_2^{-1/\alpha_2} u^{-2/\alpha_1-2/\alpha_2}.$$
We derive the uniform asymptotics, as $u\rw\IF$, of
$$\pk{ Vol(\{ t\in I_{k_1,k_2}(u,n): X(t) >u \})> \vv(u)x}=\pk{\int_{[0,n]^2} \mathbb{I} ( \xi_{u,n,k_1,k_2}(t)> u) dt>x},$$
with $x\geq 0$.
For this, we  check conditions {\bf C0-C3} of  Lemma \ref{the-weak-conv}
with $\Gamma(f)=f, f\in C([0,n]^2)$. First note that {\bf C0-C1} follow trivially with $h=0$ and $g_{u,j}=u$. Moreover, by (\ref{cor}), we have
\BQNY
\lim_{u\rw\IF}\sup_{0\leq k_i\leq N_{i}(u,n), i=1,2}\sup_{s,t\in [0,n]^2}\left|u^2Var(\xi_{u,n,k_1,k_2}(t)-\xi_{u,n,k_1,k_2}(s))-2Var\left(\sum_{i=1}^2B_{\alpha_i}(t_i)-\sum_{i=1}^2B_{\alpha_i}(s_i)\right)\right|=0,
\EQNY
with $B_{\alpha_i}, i=1,2$ being two independent fBms' with indices $\alpha_i/2$, respectively. This implies that {\bf C2} is satisfied with
$\zeta(t)=\sum_{i=1}^2B_{\alpha_i}(t_i).$
Additionally, in light of (\ref{cor}), we have that
$$\sup_{0\leq k_i\leq N_{i}(u,n)+1, i=1,2}u^2Var(\xi_{u,n,k_1,k_2}(t)-\xi_{u,n,k_1,k_2}(s))\leq C\norm{\vk{t}-\vk{s}}^{\min(\alpha_1, \alpha_2)}, \quad \vk{s},\vk{t} \in [0,n]^2.$$
This means that  {\bf C3}  holds.
Thus, by Lemma \ref{the-weak-conv},
\BQN\label{uniform1}
\lim_{u\rw\IF}\sup_{0\leq k_i\leq N_{i}(u,n), i=1,2}\left|\frac{\pk{ Vol(\{ t\in I_{k_1,k_2}(u,n): X(t) >u \})> \vv(u)x}}{\Psi(u)}-\mathcal{B}_{\alpha_1,\alpha_2}(x, [0,n]^2)\right|=0.
\EQN
Therefore, by Lemma 6.1 in \cite{Pit96} we obtain
$$\lim_{u\rw\IF}\sup_{0\leq k_i\leq N_{i}(u,n), i=1,2}\left|\frac{\pk{ Vol(\{ t\in I_{k_1,k_2}(u,n): X(t) >u \})> \vv(u)x}}{\pk{\sup_{t\in I_{k_1,k_2}(u,n)}X(t)>u}}-\frac{\mathcal{B}_{\alpha_1,\alpha_2}(x, [0,n]^2)}{\mathcal{B}_{\alpha_1,\alpha_2}(0,[0,n]^2)}\right|=0.$$

Since, by (i) of Lemma \ref{l.1}, for any $x\geq 0$ we have
\BQN \label{ehe}
\mathcal{B}_{\alpha_1,\alpha_2}(x)=\lim_{n\rw\IF}\frac{\mathcal{B}_{\alpha_1,\alpha_2}(x, [0,n]^2)}{n^2}\in (0,\IF),
\EQN
then
\BQN\label{Constant}
\frac{\mathcal{B}_{\alpha_1,\alpha_2}(x)}{\mathcal{B}_{\alpha_1,\alpha_2}(0)}=\lim_{n\rw\IF}\frac{\mathcal{B}_{\alpha_1,\alpha_2}(x, [0,n]^2)}{\mathcal{B}_{\alpha_1,\alpha_2}(0,[0,n]^2)}\in (0,{1]},  \quad x\geq 0,
\EQN
which confirms that A2 holds with
$\bar F(x)=\frac{\mathcal{B}_{\alpha_1,\alpha_2}(x)}{\mathcal{B}_{\alpha_1,\alpha_2}(0)}$.

{\it \underline{Condition A3.}}
By (7.4) in the proof of Lemma 7.1 in \cite{Pit96},  for all large $u$ and $n$
\begin{eqnarray*}
\sum_{0\leq k_i, k_i'\leq N_i(u,n), i=1,2, (k_1,k_2)\neq (k_1',k_2')}\pk{\sup_{t\in I_{k_1,k_2}(u,n)}X(t)>u, \sup_{t\in I_{k_1',k_2'}(u,n)}X(t)>u }
\leq
\left(\frac{\CC_2}{\sqrt{n}}+ e^{-\CC_1n^{\CC}}\right)\pk{ \sup_{t\in E} X(t)>u},
\end{eqnarray*}
where $\CC,\CC_1$ and $\CC_2$ are some positive constants,
which gives that A3 is satisfied.

This completes the proof. \QED

\subsection{Proof of Proposition \ref{Onepoint}}
 Without loss of generality, we assume that $t^*=(0,0)$. 
 The proof relies on verification that A1-A3 in Theorem \ref{th1} are satisfied.
 We begin by introducing  some notation. Let \BQN\label{Ik}
I_{k_1,k_2}(u,n)=\prod_{i=1}^2[k_iv_i(u)n, (k_i+1)v_i(u)n],  \quad v_i(u)=a_i^{-1/\alpha_i^*} u^{-2/\min(\alpha_i,\beta_i)}, i=1,2,\quad v(u)=v_1(u)v_2(u),
\EQN
where $\alpha_i^*=\alpha_i \I(\alpha_i\leq \beta_i)+\IF \I(\alpha_i> \beta_i)$.
Additionally, let $$e(t)=\frac{1-{\sigma}(t)}{\sum_{i=1}^2b_i|t_i|^{\beta_i}}-1, |t|\neq 0, ~e_u=\sup_{0<|t_i|<\left(\frac{\ln u}{u}\right)^{2/\beta_i}} |e(t)|,$$
and set $$ N_{i}'(u,n)=
\left[\frac{(e_u^{-1/4}\wedge\ln u)^{2/\beta_i}}{u^{2/\beta_i}v_i(u)n}\right],\ i=1,2. $$

We
distinguish different scenarios according to the values of $\alpha_i, \beta_i, i=1,2$.
\\
{\bf \underline{Case $\alpha_i<\beta_i, i=1,2$}}. In this scenario
$$v_i(u)=a_i^{-1/\alpha_i} u^{-2/\alpha_i}, i=1,2, \ \
K_{u,n}=\{ (k_1,k_2): 0\leq |k_i|\leq N_{i}'(u,n),i=1,2\}$$
and
$E(u,n)=\bigcup_{(k_1,k_2)\in K_{u,n}}I_{k_1,k_2}(u,n)$.

{\it \underline{Conditions A1 and A3.}}
Following the same reasoning as in the proof of Proposition \ref{THSTA},
the validity of conditions A1 and A3 follows straightforwardly from
(34), (40) and (41) in \cite{KEP20151}.

{\it \underline{Condition A2.}}
Let $$\xi_{u,n,k_1,k_2}(t)=\overline{X}(v_1(u)(k_1n+t_1), v_2(u)(k_2n+t_2)),$$
\BQN\label{uk}
u_{n,k_1,k_2}^{-}=u\inf_{t\in I_{k_1,k_2}(u,n)}\frac{1}{{\sigma(t)}},
\quad u_{n,k_1,k_2}^{+}=u\sup_{t\in I_{k_1,k_2}(u,n)}\frac{1}{{\sigma(t)}}.
\EQN
Then
\BQNY\pk{ Vol(\{ t\in I_{k_1,k_2}(u,n):X(t) >u \})\ge  \vv(u)x}\leq \pk{\int_{[0,n]^2} \mathbb{I} ( \xi_{u,n,k_1,k_2}(t)> u_{n,k_1,k_2}^{-}) dt>x},\\
\pk{ Vol(\{ t\in I_{k_1,k_2}(u,n):X(t) >u \})\ge  \vv(u)x}\geq \pk{\int_{[0,n]^2} \mathbb{I} (\xi_{u,n,k_1,k_2}(t)> u_{n,k_1,k_2}^{+}) dt>x}.
\EQNY
In order to derive the uniform asymptotics of the above terms
we  check conditions {\bf C0-C3} of  Lemma \ref{the-weak-conv}
with $\Gamma(f)=f,$  $f\in C([0,n]^2)$
for $\xi_{u,n,k_1,k_2}(t), \ (k_1,k_2)\in K_{u,n}.$

Note that {\bf C0-C1} holds with $h=0$ and $g_{u,j}=u_{n,k_1,k_2}^{\pm}$.  By (\ref{var}) and (\ref{cor2}), we have
\BQNY
\lim_{u\rw\IF}\sup_{s,t\in [0,n]^2, (k_1,k_2)\in K_{u,n} }\left|(u_{n,k_1,k_2}^{\pm})^2(Var(\xi_{u,n,k_1,k_2}(t)-\xi_{u,n,k_1,k_2}(s)))-2Var\left(\sum_{i=1}^2B_{\alpha_i}
(t_i)-\sum_{i=1}^2B_{\alpha_i}(s_i)\right)\right|=0,
\EQNY
where $B_{\alpha_i}, i=1,2$ are two independent fBm's with indices $\alpha_i, i=1,2$ respectively.
This confirms that {\bf C2} holds with $\zeta(t_1,t_2)=B_{\alpha_1}(t_1)+B_{\alpha_2}(t_2)$.  By (\ref{cor2}), we have
$$\sup_{(k_1,k_2)\in K_{u,n}}(u_{n,k_1,k_2}^{\pm})^2(Var(\xi_{u,n,k_1,k_2}(t)-\xi_{u,n,k_1,k_2}(s)))\leq Q||s-t||^{\min(\alpha_1,\alpha_2)}, \quad s,t\in [0,n]^2.$$
Thus {\bf C3} is satisfied.

Therefore, by Lemma \ref{the-weak-conv}, we have that for $0\leq x<n^2$,
\BQN\label{uniform3}\lim_{u\rw\IF}\sup_{(k_1,k_2)\in K_{u,n}}\left|\frac{\pk{\int_{[0,n]^2} \mathbb{I} ( \xi_{u,n,k_1,k_2}(t)> u_{n,k_1,k_2}^{\pm}) dt>x}}{\Psi(u_{n,k_1,k_2}^{\pm})}-\mathcal{B}_{\alpha_1,\alpha_2}(x, [0,n]^2)\right|=0.
\EQN
Since
\BQN\label{upper-lower}
\lim_{u\rw\IF}\sup_{(k_1,k_2)\in K_{u,n}}\left|\frac{\Psi(u_{n,k_1,k_2}^{-})}{\Psi(u_{n,k_1,k_2}^{+})}-1\right|=0
\EQN
(see Section \ref{s.supl} for the validation of (\ref{upper-lower})),
{by} (\ref{uniform3}) we obtain
for $0\leq x<n^2$
$$\lim_{u\rw\IF}\sup_{(k_1,k_2)\in K_{u,n}}\left|\frac{\pk{ Vol(\{ t\in I_{k_1,k_2}(u,n):X(t) >u \})\ge  \vv(u)x}}{\Psi(u_{n,k_1,k_2}^{-})}-\mathcal{B}_{\alpha_1,\alpha_2}(x, [0,n]^2)\right|=0.$$
Therefore, (\ref{Pickands}) holds with
$$ \bar F_n(x)=\frac{\mathcal{B}_{\alpha_1,\alpha_2}(x,[0,n]^2)}{\mathcal{B}_{\alpha_1,\alpha_2}(0,[0,n]^2)}, \quad x\geq 0.$$
Finally, by (\ref{Constant}), we have that A2 holds. Thus the claim is established with
$$ \bar F(x)=\frac{\mathcal{B}_{\alpha_1,\alpha_2}(x)}
 {\mathcal{B}_{\alpha_1,\alpha_2}(0)}.$$

 {\bf \underline{Case $\alpha_1=\beta_1, \alpha_2<\beta_2$}}. In this case
 $ v_i(u)=a_i^{-1/\alpha_i} u^{-2/\alpha_i}, i=1,2.$
Let
\BQN\label{E1}
 \hat{I}_{k_2}(u,n)=I_{-1,k_2}(u,n)\cup I_{0,k_2}(u,n),\quad  E_1(u,n)=\bigcup_{k_2\in K_{u,n}}\hat{I}_{k_2}(u,n),
 \EQN
 where
 $K_{u,n}\coloneqq\{ k_2\in \mathbb{Z}:|k_2|\leq N_2'(u,n)\}$.

{\it \underline{Conditions A1 and A3.}}
Analogously to the previous case, conditions A1 and A3 hold with
$E(u,n)\coloneqq E_1(u,n)$ and ${I}_{k}(u,n)\coloneqq \hat{I}_{k_2}(u,n)$,
by
(34), (46), (48) and (49) of \cite{KEP20151}.

{\it \underline{Condition A2.}}
 Rewrite (\ref{var}) as
$$
\frac{1}{\sigma(t)}=\left(1+(1+e_1(t_1))b_1|t_1|^{\beta_1}\right)\left(1+(1+e_2(t_2))b_2|t_2|^{\beta_2}\right),
$$
{for some functions $e_1(t_1)$ and $e_2(t_2)$ which satisfy}
 $$\lim_{u\rw\IF}\sup_{t\in E_1(u,n)}|e_i(t_i)|=0,\quad  i=1,2.$$
Let  $$\xi_{u,n,k_2}(t)=\frac{\overline{X}(v_1(u)t_1, v_2(u)(k_2n+t_2))}{1+b_1|v_1(u)t_1|^{\beta_1}(1+e_1(v_1(u)t_1))}, \quad  v(u)=a_1^{-1/\alpha_1}a_2^{-1/\alpha_2} u^{-2/\alpha_1-2/\alpha_2},$$
 $$u_{k_2,n}^{-}=u\inf_{t\in \hat{I}_{k_2}(u,n)}(1+b_2|t_2|^{\beta_2}(1+e_2(t_2))),\quad u_{k_2,n}^{+}=u\sup_{t\in \hat{I}_{k_2}(u,n)}(1+b_2|t_2|^{\beta_2}(1+e_2(t_2))).$$
Then it follows that
\BQNY\pk{ Vol(\{ t\in \hat{I}_{k_2}(u,n):X(t) >u \})> \vv(u)x}\leq \pk{\int_{[-n,n]\times[0,n]} \mathbb{I} ( \xi_{u,n,k_2}(t)> u_{k_2,n}^{-}) dt>x},\\
\pk{ Vol(\{ t\in \hat{I}_{k_2}(u,n):X(t) >u \})>  \vv(u)x}\geq \pk{\int_{[-n,n]\times[0,n]} \mathbb{I} ( \xi_{u,n,k_2}(t)> u_{k_2,n}^{+}) dt>x}.
\EQNY
Straightforward application of
Lemma \ref{the-weak-conv} with $\Gamma(f)=f,$  $f\in C([-n,n]\times [0,n])$
and
$h(t)=a_1^{-1}b_1|t_1|^{\alpha_1}$ in {\bf C1},
gives that
%
for $0\leq x<2n^2$,
 \BQNY
\lim_{u\rw\IF}\sup_{ k_2\in K_{u,n}}\left|\frac{\pk{\int_{[-n,n]\times[0,n]} \mathbb{I} ( \xi_{u,n,k_2}(t)> u_{k_2,n}^{\pm}) dt>x
}}{\Psi(u_{k_2,n}^{\pm})}-\mathcal{B}_{\alpha_1,\alpha_2}^{a_1^{-1}b_1|t_1|^{\alpha_1},0}(x, [-n,n]\times[0,n])\right|=0.
\EQNY
Similarly to (\ref{upper-lower}), we have
\BQN\label{upper-lower1}
\lim_{u\rw\IF}\sup_{k_2\in K_{u,n}}\left|\frac{\Psi(u_{k_2,n}^{-})}{\Psi(u_{k_2,n}^{+})}-1\right|=0.
\EQN
Consequently, for $0\leq x<2n^2$
\BQN\label{uniform2}\lim_{u\rw\IF}\sup_{k_2\in K_{u,n}}\left|\frac{\pk{ Vol(\{ t\in \hat{I}_{k_2}(u,n):X(t) >u \})>  \vv(u)x}}{\Psi(u_{k_2,n}^{-})}-\mathcal{B}_{\alpha_1,\alpha_2}^{a_1^{-1}b_1|t_1|^{\alpha_1},0}(x, [-n,n]\times[0,n])\right|=0.\nonumber\\
\EQN
Thus (\ref{Pickands}) holds with
$$\bar F_n(x)=\frac{\mathcal{B}_{\alpha_1,\alpha_2}^{a_1^{-1}b_1|t_1|^{\alpha_1},0}(x, [-n,n]\times[0,n])}{\mathcal{B}_{\alpha_1,\alpha_2}^{a_1^{-1}b_1|t_1|^{\alpha_1},0}(0, [-n,n]\times[0,n])}.$$
By (ii) of Lemma \ref{l.1} it follows that
\begin{equation}\label{limit} \lim_{n\to\infty}\bar F_n(x)=\frac{\mathcal{B}_{\alpha_1,\alpha_2}^{a_1^{-1}b_1|t_1|^{\alpha_1},0}(x)}{\mathcal{B}_{\alpha_1,\alpha_2}^{a_1^{-1}b_1|t_1|^{\alpha_1},0}(0) }\in (0,{1]},
\end{equation}
which confirms that A2 holds.
Thus, applying Theorem \ref{th1}, we establish the claim with
$$ \bar F(x)=\frac{\mathcal{B}_{\alpha_1,\alpha_2}^{a_1^{-1}b_1|t_1|^{\alpha_1},0}(x)}{\mathcal{B}_{\alpha_1,\alpha_2}^{a_1^{-1}b_1|t_1|^{\alpha_1},0}(0) }.$$
 {\bf \underline{Case $\alpha_1=\beta_1, \alpha_2=\beta_2$}}. In this case we have
 $ v_i(u)=a_i^{-1/\alpha_i} u^{-2/\alpha_i}, i=1,2.$
 Let
 \BQN\label{EE1}
 E(u,n)\coloneqq
 \hat{I}(u,n)\coloneqq\bigcup_{i,j=-1,0}I_{i,j}(u,n).
 \EQN
{\it \underline{Conditions A1 and A3.}}
 It follows
from (34) and (52) in the proof of  theorem 3.1 of \cite{KEP20151} that
A1 holds. Since {we take}
only one interval {$I_1(u,n)$}, condition A3 is not applicable to this case.

{\it \underline{Condition A2.}}
 \COM{(\ref{var}) can be rewritten as
 $$\sigma(t)=1+b_1|t_1|^{\beta_1}(1+e_3(t))+b_2|t_2|^{\beta_2}(1+e_4(t)),$$
 where for any $n>0$
 $$\lim_{u\rw\IF}\sup_{t\in \hat{I}(u,n)}|e_i(t)|=0, \quad i=3,4.$$}
 Let
 \BQNY
 \xi_{u,n}(t)=X(v_1(u)t_1, v_2(u)t_2)), \quad v(u)=a_1^{-1/\alpha_1}a_2^{-1/\alpha_2} u^{-2/\alpha_1}u^{-2/\alpha_2}.
 \EQNY
 Then
 $$\pk{ Vol(\{ t\in \hat{I}(u,n): X(t) >u \})\ge  \vv(u)x}=\pk{\int_{[-n,n]^2} \mathbb{I} ( \xi_{u,n}(t)> u) dt>x}.$$
 In order to derive the  asymptotics of the above term,
similarly to the previous cases, we observe that
{\bf C1} in Lemma \ref{the-weak-conv} holds with
$h(t)=a_1^{-1}b_1|t_1|^{\alpha_1}+a_2^{-1}b_2|t_2|^{\alpha_2}$
while
{\bf C2} and {\bf C3} have been checked in the case of $\alpha_i<\beta_i, i=1,2$.

Hence we have
 \BQNY
 \lim_{u\rw\IF}\left|\frac{\pk{\int_{[-n,n]^2} \mathbb{I} ( \xi_{u,n}(t)> u) dt>x
}}{\Psi(u)}-\mathcal{B}_{\alpha_1,\alpha_2}^{a_1^{-1}b_1|t_1|^{\alpha_1},a_2^{-1}b_2|t_2|^{\alpha_2}}(x, [-n,n]^2)\right|=0.
 \EQNY
Combining the above with the fact that, by (iii) of Lemma \ref{l.1},
 \BQN\lim_{n\rw\IF} \mathcal{B}_{\alpha_1,\alpha_2}^{a_1^{-1}b_1|t_1|^{\alpha_1},a_2^{-1}b_2|t_2|^{\alpha_2}}(x, [-n,n]^2)\in (0,\IF)
\nonumber
 \EQN
we conclude that A2 holds with
 $$ \bar F(x)=\frac{ \mathcal{B}_{\alpha_1,\alpha_2}^{a_1^{-1}b_1|t_1|^{\alpha_1},a_2^{-1}b_2|t_2|^{\alpha_2}}(x)}{ \mathcal{B}_{\alpha_1,\alpha_2}^{a_1^{-1}b_1|t_1|^{\alpha_1},a_2^{-1}b_2|t_2|^{\alpha_2}}(0)}\in (0,{1]}.$$
 Hence we establish the claim. \\

 For the cases $\alpha_1>\beta_1, \alpha_2=\beta_2$, and  $\alpha_1>\beta_1, \alpha_2>\beta_2$, we can establish the claim similarly to the  case of $\alpha_1=\beta_1, \alpha_2=\beta_2$. For the case $\alpha_1>\beta_1, \alpha_2<\beta_2$, the proof is similar to the case of $\alpha_1=\beta_1$, $\alpha_2<\beta_2$.
 This completes the proof.\QED

 \subsection{Proof of Proposition \ref{chith}}
 In order to apply Theorem  \ref{th1}, we introduce some useful notation. Let
\BQNY
I_{k}(u,n)=[kv(u)n, (k+1)v(u)n], \quad N(u,n)=\left[\frac{T}{v(u)n}\right]-1,
\EQNY
and
$E(u,n)=\bigcup_{k\in K_{u,n}} I_{k}(u,n),$ with
$K_{u,n}=\{ k\in \mathbb{N}: 0\leq k\leq N(u,n)  \}$ and $ v(u)=a^{-1/\alpha}u^{-2/\alpha}.$
We denote by
$$Z(t,\theta)=\sum_{i=1}^m X_i(t)v_i(\theta), \quad  A=[0,\pi]^{m-2}\times[0,2\pi),$$
where $\theta=(\theta_1, \ldots, \theta_{m-1})$ and
$$v_1(\theta)=\cos\theta_1, v_2(\theta)=\sin\theta_1\cos\theta_2, v_3(\theta)=\sin\theta_1\sin\theta_2\cos\theta_3, \dots, v_{m-1}(\theta)=(\prod_{i=1}^{m-2}\sin\theta_i)\cos\theta_{m-1}, v_{m}(\theta)=\prod_{i=1}^{m-1}\sin\theta_i.$$
In this proof, we will use that
$$\chi(t)
=\sup_{\theta\in A} Z(t,\theta).$$
We split the set $A$ into (setting $k=(k_1,\dots, k_{m-1})$)
$$A=\bigcup_{k\in \Lambda}A_{k}, \quad \Lambda=\{(k_1,\dots, k_{m-1}): 1\leq k_i \leq L, 1\leq i\leq m-2, 1\leq k_{m-1}\leq 2L\},$$
where
$$A_{k}=\prod_{i=1}^{m-1}\left[\frac{(k_i-1)\pi}{L}, \frac{k_i\pi}{L}\right], \quad k_{m-1}\leq 2L-1,$$
$$
A_{k_1,\dots, k_{m-2}, 2L}=\left(\prod_{i=1}^{m-2}\left[\frac{(k_i-1)\pi}{L}, \frac{k_i\pi}{L}\right]\right)\times \left[2\pi-\frac{\pi}{L}, 2\pi\right),$$
and $L$ is a positive integer. 
Moreover,  let
\BQN\label{pi5}
\pi_1(u)&\coloneqq&\sum_{k\neq k', k,k'\in \Lambda}\pk{\sup_{t\in [0,v(u)n], \theta\in A_{k}} Z(t,\theta)>u, \sup_{t\in [0,v(u)n], \theta\in A_{k'}} Z(t,\theta)>u},\\
 \Sigma\Sigma_{u,n}
 &\coloneqq&\sum_{0\leq k_1<k_2\leq N(u,n)}\pk{\sup_{t\in I_{k_1}(u,n)}\chi(t)>u, \sup_{t\in I_{k_2}(u,n)}\chi(t)>u}\nonumber\\
 &=& \sum_{0\leq k_1<k_2\leq N(u,n)}\pk{\sup_{(t,\theta)\in I_{k_1}(u,n)\times A}Z(t,\theta)>u, \sup_{(t,\theta)\in I_{k_2}(u,n)\times A}Z(t,\theta)>u}\nonumber\\
 &\leq&\sum_{0\leq k_1<k_2\leq N(u,n), i,j\in \Lambda}\pk{\sup_{(t,\theta)\in I_{k_1}(u,n)\times A_i}Z(t,\theta)>u, \sup_{(t,\theta)\in I_{k_2}(u,n)\times A_j}Z(t,\theta)>u}.\nonumber
\EQN
Denote by (with $k=(k_1,\dots, k_{m-1}), l=(l_1,\dots,l_{m-1})$)
$$J_{k,l}(u)=\prod_{i=1}^{m-1}\left[\frac{(k_i-1)\pi}{L}+l_iu^{-1}n_1, \frac{(k_i-1)\pi}{L}+(l_i+1)u^{-1}n_1\right], \quad \Lambda_1(u)=\left\{l: 0\leq l_i\leq \left[\frac{\pi u}{Ln_1}\right], 1\leq i\leq m-1\right\},$$
and let
\BQN\label{pk}
p_k^*(u)=\sum_{l,l'\in\Lambda_1(u), l\neq l'}\pk{\sup_{t\in [0,v(u)n],\theta\in J_{k,l}(u)}Z(t,\theta)>u, \sup_{t\in [0,v(u)n],\theta\in J_{k,l'}(u)}Z(t,\theta)>u }.
\EQN

{\it \underline{Conditions A1 and A3.}}
Condition A1   follows from
 Corollary 7.3 in \cite{Pit96} while A3  can be deduced from
equations (7.4), (7.6) and (7.18) in the proofs of Lemma 7.1 and Theorem 7.1 in \cite{Pit96}.

{\it \underline{Condition A2.}}
Let us put
\BQNY
\pi(n,u)\coloneqq \pk{\int_{[0,v(u)n]}\I(\chi(t)>u)dt>\vv(u)x}=\pk{\int_{[0,v(u)n]}\I\left(\sup_{\theta\in A}Z(t,\theta)>u\right)dt>\vv(u)x}.
\EQNY
To verify A2, {by stationarity} we have to find the asymptotics of $\pi{(n,u)}$ as $u\to\infty$,
which is given in the following lemma.
\BEL\label{l.xi}
For $n>x$
\BQNY
\pi(n,u)\sim \frac{\widehat{\mathcal{B}}_{\alpha,2,\dots,2}(x, n)}{\widehat{\mathcal{B}}_{\alpha,2,\dots,2}(0, n)}\pk{\sup_{[0, v(u)n]}\chi(t)>u},\quad u\rw\IF.
\EQNY
\EEL
{\it Proof of Lemma \ref{l.xi}}. Let
$D_k=\{t\in [0, v(u)n]: \sup_{\theta\in A_k}Z(t,\theta)>u \}.$
Then we have
\BQNY
\int_{[0,v(u)n]}\I\left(\sup_{\theta\in A}Z(t,\theta)>u\right)dt&=&\int_{[0,v(u)n]}\I_{\bigcup_{k\in\Lambda} D_k}(t)dt\\
&\leq&\sum_{k\in\Lambda}\int_{[0,v(u)n]}\I_{D_k}(t)dt
\EQNY
and
\BQNY
\int_{[0,v(u)n]}\I\left(\sup_{\theta\in A}Z(t,\theta)>u\right)dt&\geq&\sum_{k\in\Lambda}\int_{[0,v(u)n]}\I_{D_k}(t)dt-\sum_{k\neq k', k, k'\in \Lambda} \int_{[0,v(u)n]}\I_{D_k\bigcap D_{k'}}(t)dt.
\EQNY
Note that
\BQNY
\pi(n,u)&\geq& \mathbb{P}\left(\sum_{k\in\Lambda}\int_{[0,v(u)n]}\I_{D_k}(t)dt-\sum_{k\neq k', k, k'\in \Lambda} \int_{[0,v(u)n]}\I_{D_k\bigcap D_{k'}}(t)dt>v(u)x\right)\\
&\geq& \mathbb{P}\left(\sum_{k\in\Lambda}\int_{[0,v(u)n]}\I_{D_k}(t)dt>v(u)(x+\epsilon), \sum_{k\neq k', k, k'\in \Lambda} \int_{[0,v(u)n]}\I_{D_k\bigcap D_{k'}}(t)dt\leq v(u)\epsilon\right)\\
&\geq& \mathbb{P}\left(\sum_{k\in\Lambda}\int_{[0,v(u)n]}\I_{D_k}(t)dt>v(u)(x+\epsilon)\right)-\mathbb{P}\left(\sum_{k\neq k', k, k'\in \Lambda} \int_{[0,v(u)n]}\I_{D_k\bigcap D_{k'}}(t)dt>v(u)\epsilon\right)\\
&\geq&\mathbb{P}\left(\sum_{k\in\Lambda}\int_{[0,v(u)n]}\I_{D_k}(t)dt>v(u)(x+\epsilon)\right)-\pi_1(u)\\
&\geq& \sum_{k\in \Lambda^*}p_{k}(x+\epsilon,u)-2\pi_1(u),
\EQNY
where $\epsilon>0$ and $\pi_1(u)$ is given in (\ref{pi5}) and
\BQNY
p_{k}(x,u)&=&\pk{\int_{[0,v(u)n]}\I\left(\sup_{\theta\in A_k}Z(t,\theta)>u\right)dt>\vv(u)x},\\
 \Lambda^*&=&\{k\in \Lambda, 1<k_i<L, 1\leq i\leq m-2, k_{m-1}\neq 1, L, 2L\}.
\EQNY
Similarly we get
$$\pi(n,u)\leq \sum_{k\in \Lambda}p_{k}(x,u)+\pi_1(u).$$
Hence
\BQN\label{upper-lower2}
\sum_{k\in \Lambda^*}p_{k}(x+\epsilon,u)-2\pi_1(u)\leq \pi(n,u)\leq \sum_{k\in \Lambda}p_{k}(x,u)+\pi_1(u).
\EQN
{$\diamond$ \it \underline{Upper bound for $p_{k}(x,u)$}}. A direct calculations show
\BQNY
Var(Z(t,\theta))&=&1,\\
 Corr(Z(t,\theta), Z(t',\theta'))&=&Corr(X(t), X(t'))\left(\cos(\theta_1-\theta_1')-\sin\theta_1\sin\theta _1'(1-\cos(\theta_2-\theta_2'))\right.\\
 &&\quad \left.-\dots-
\left(\prod_{i=1}^{m-2}\sin\theta_i\sin\theta_i'\right)\left(1-\cos(\theta_{m-1}-\theta_{m-1}')\right)\right).
\EQNY
Hence
\BQN\label{corchi}
1-Corr(Z(t,\theta), Z(t',\theta'))
 &\sim& a|t-t'|^{\alpha}+\frac{1}{2}(\theta_1-\theta_1')^2+\frac{\sin^2\theta_1}{2}(\theta_2-\theta_2')^2\nonumber\\
 &&\quad +\frac{1}{2}
\left(\prod_{i=1}^{m-2}\sin^2\theta_i\right)(\theta_{m-1}-\theta_{m-1}')^2, \quad |t-t'|\rw 0, ||\theta-\theta'||\rw 0.
\EQN
We have
\BQN\label{upper5}
p_{k}(x,u)\leq \sum_{l\in\Lambda_1(u)}\pk{\int_{[0,v(u)n]}\I\left(\sup_{\theta\in J_{k,l}(u)}Z(t,\theta)>u\right)dt>\vv(u)x}+p_k^*(u),
\EQN
where $p_k^*(u)$ is given in (\ref{pk}).
Let $$Z_{u,k,l}(t,\theta)=Z\left(v(u)t, \frac{(k_1-1)\pi}{L}+l_1u^{-1}n_1+u^{-1}c_1^{-1}(\theta_{k,l}(u))\theta_1, \dots,  \frac{(k_{m-1}-1)\pi}{L}+l_{m-1}u^{-1}n_{1}+u^{-1}c_{m-1}^{-1}(\theta_{k,l}(u))\theta_{m-1}\right),$$
and
$G_{l}=\prod_{i=1}^{m-1}[0,c_i(\theta_{k,l}(u))n_1],$
where $$c_k(\theta)=2^{-1/2}\prod_{i=1}^{k-1}|\sin\theta_i|, 2\leq k\leq m-1,\quad  c_1(\theta)=2^{-1/2}, \quad \theta_{k,l}(u)=\left(\frac{(k_1-1)\pi}{L}+l_1u^{-1}n_1,\dots,\frac{(k_{m-1}-1)\pi}{L}+l_{m-1}u^{-1}n_1\right).$$
Noting that
$$G_{l}=\prod_{i=1}^{m-1}[0,c_i(\theta_l(u))n_1]\subset\prod_{i=1}^{m-1}[0,c_{k,i}^+n_1]=:G_k^+, \quad c_{k,i}^+=\sup_{\theta\in A_{k}}c_i(\theta), $$
we have
\BQNY
\pk{\int_{[0,v(u)n]}\I\left(\sup_{\theta\in J_{k,l}(u)}Z(t,\theta)>u\right)dt>\vv(u)x}&=&\pk{\int_{[0,n]}\I\left(\sup_{\theta\in G_{l}}Z_{u,k,l}(t,\theta)>u\right)dt>\vv(u)x}\\
&\leq&\pk{\int_{[0,n]}\I\left(\sup_{\theta\in G_{k}^+}Z_{u,k,l}(t,\theta)>u\right)dt>\vv(u)x}.
\EQNY

A straightforward application of Lemma \ref{the-weak-conv}
for
$\Gamma: C([0,n]\times G_k^+)\rw C([0,n])$ defined by
$\Gamma(f)=\sup_{\theta\in G_k^+}f(t,\theta), \quad f\in C([0,n]\times G_k^+)$, where
$h=0$ in {\bf C1} and
$\zeta(t,\theta)=B_{\alpha}(t)+\sum_{i=1}^{m-1}N_i\theta_i,$
with $N_i, i=1,\dots, m-1$ being independent standard normal random variables independent of $B_\alpha$,
implies that
for all $x\geq 0$ we have
\BQN\label{uniformchi}
\lim_{u\rw\IF}\sup_{l\in \Lambda_1(u)}\left|\frac{\pk{\int_{[0,n]}\I\left(\sup_{\theta\in G_{l}^+}Z_{u,k,l}(t,\theta)>u\right)dt>\vv(u)x}}{\Psi(u)}-
\widehat{\mathcal{B}}_{\alpha,2,\dots, 2}(x, [0,n]\times G_l^+)\right|=0.
\EQN
By (7.18) in the proof of  Theorem 7.1 in \cite{Pit96}, we have
\BQN\label{upperchi}
p_k^*(u)= o\left(u^{m-1}\Psi(u)\right), \quad u\rw\IF, n_1\rw\IF.
\EQN
Hence, by (\ref{upper5})-(\ref{upperchi}) and using  Lemma \ref{lemA0} we have
\BQNY\label{upperchi1}
p_{k}(x,u)&\leq& \limsup_{n_1\rw\IF}\frac{\widehat{\mathcal{B}}_{\alpha,2,\dots,2}(x, [0,n]\times G_l^+)}{(n_1)^{m-1}}
\left(\frac{\pi}{L}\right)^{m-1}u^{m-1}\Psi(u)\\
&\leq &\widehat{\mathcal{B}}_{\alpha,2,\dots,2}(x, n)\prod_{i=1}^{m-1}c_{k,i}^+
\left(\frac{\pi}{L}\right)^{m-1}u^{m-1}\Psi(u), \quad u\rw\IF, n_1\rw\IF.
\EQNY
{$\diamond$ \it \underline{Lower bound for $p_{k}(x,u)$}}.
By (\ref{upper-lower2}), we have  that for $\epsilon>0$
\BQNY
p_{k}(x,u)
 &\geq&\sum_{l\in\Lambda_2(u)}\pk{\int_{[0,v(u)n]}\I\left(\sup_{\theta\in J_{{k,l}}(u)}Z(t,\theta)>u\right)dt>\vv(u)(x+\epsilon)}-2p_k^*(u)\\
 &\geq&\sum_{l\in\Lambda_2(u)}\pk{\int_{[0,n]}\I\left(\sup_{\theta\in G_{k}^-}Z_{u,k,l}(t,\theta)>u\right)dt>\vv(u)(x+\epsilon)}-2p_k^*(u),
\EQNY
where $$\Lambda_2(u)=\left\{l: 0\leq l_i\leq \left[\frac{\pi u}{Ln_1}\right]-1, 1\leq i\leq m-1\right\},$$ $$G_{l}=\prod_{i=1}^{m-1}[0,c_i(\theta_l(u))n_1]\supset\prod_{i=1}^{m-1}[0,c_{k,i}^-n_1]=:G_k^-, \quad c_{k,i}^-=\min_{\theta\in A_{k}}c_i(\theta). $$
By (\ref{uniformchi}), (\ref{upperchi}), Lemma \ref{lemA0} and Remark \ref{lemma} we have for $n>x$
\BQNY\label{lowerchi}
p_{k}(x,u)&\geq& \liminf_{n_1\rw\IF}\frac{\widehat{\mathcal{B}}_{\alpha,2,\dots,2}(x+\epsilon, [0,n]\times G_l^-)}{(n_1)^{m-1}}
\left(\frac{\pi}{L}\right)^{m-1}u^{m-1}\Psi(u)\\
&\geq&\widehat{\mathcal{B}}_{\alpha,2,\dots,2}(x+\epsilon, n)\prod_{i=1}^{m-1}c_{k,i}^-
\left(\frac{\pi}{L}\right)^{m-1}u^{m-1}\Psi(u)\\
&\geq& \widehat{\mathcal{B}}_{\alpha,2,\dots,2}(x, n)\prod_{i=1}^{m-1}c_{k,i}^-
\left(\frac{\pi}{L}\right)^{m-1}u^{m-1}\Psi(u), \quad u\rw\IF, 
\epsilon\rw 0.
\EQNY

{$\diamond$ \it \underline{Asymptotics for $\pi(u,n)$}}.
By (7.6) in \cite{Pit96} 
$$\pi_1(u)=o\left(u^{m-1}\Psi(u)\right), \quad u\rw\IF, L\rw\IF.$$
Therefore, in view of (\ref{upper-lower2}),
\BQNY
 \limsup_{u\to\infty}\frac{\pi(n,u)}{u^{m-1}\Psi(u)}&\leq& \limsup_{L\rw\IF}\sum_{k\in\Lambda}\left(\prod_{i=1}^{m-1}c_{k,i}^+\right)\left(\frac{\pi}{L}\right)^{m-1} \widehat{\mathcal{B}}_{\alpha,2,\dots,2}(x,n),\\
 \liminf_{u\to\infty}\frac{\pi(n,u)}{u^{m-1}\Psi(u)}&\geq&\liminf_{L\rw\IF}\sum_{k\in\Lambda^*}\left(\prod_{i=1}^{m-1}c_{k,i}^-\right)\left(\frac{\pi}{L}\right)^{m-1} \widehat{\mathcal{B}}_{\alpha,2,\dots,2}(x, n).
\EQNY
Using the fact that
$$\limsup_{L\rw\IF}\sum_{k\in\Lambda}\left(\prod_{i=1}^{m-1}c_{k,i}^+\right)\left(\frac{\pi}{L}\right)^{m-1} =\liminf_{L\rw\IF}\sum_{k\in\Lambda^*}\left(\prod_{i=1}^{m-1}c_{k,i}^-\right)\left(\frac{\pi}{L}\right)^{m-1}=Vol(S_{m-1}), $$
it follows that
\BQNY
\pi(n,u)\sim \frac{\widehat{\mathcal{B}}_{\alpha,2,\dots,2}(x, n)}{\widehat{\mathcal{B}}_{\alpha,2,\dots,2}(0, n)}\pk{\sup_{[0, v(u)n]}\chi(t)>u},\quad u\rw\IF.
\EQNY This completes the proof of Lemma \ref{l.xi}.
\QED
{\it \underline{Condition A2 continued}}.
Lemma \ref{lemA0} yields that for $x\geq 0$
$$\frac{\widehat{\mathcal{B}}_{\alpha,2,\dots,2}(x)}{\widehat{\mathcal{B}}_{\alpha,2,\dots,2}(0)}=\lim_{n\rw\IF}\frac{\widehat{\mathcal{B}}_{\alpha,2,\dots,2}(x, [0,n])}{\widehat{\mathcal{B}}_{\alpha,2,\dots,2}(0, [0,n])}\in (0,{1]}.$$
Hence  A2 holds with
$$ \bar F(x)=\frac{\widehat{\mathcal{B}}_{\alpha,2,\dots,2}(x)}{\widehat{\mathcal{B}}_{\alpha,2,\dots,2}(0)}, \quad x\geq 0.$$ Thus we establish the claim and hence the proof is complete.\QED

\subsection{Proof of Proposition \ref{queue}}
We first apply Theorem \ref{th1} to derive the asymptotics for case ii) of Proposition \ref{queue}.
 Let $$E(u,n)=\bigcup_{i=0}^{N(u,n)} I_{i}(u,n), \quad I_i(u,n)=[iv(u)n, (i+1)v(u)n],~ N(u,n)=\left[\frac{T_u}{nv(u)}\right]-2,$$
and  $$\vv(u)=u^{\frac{2(\alpha-1)}{\alpha}}\left(\frac{(\tau^*)^{\alpha/2}}{1+c\tau^*}\right)^{2/\alpha}, \quad \tau^*=\frac{\alpha}{c(2-\alpha)}.$$
Let $$Z(s,t)=\frac{B_{\alpha}(s)-B_{\alpha}(t)}{1+c(s-t)},~I_i'(u,n)=[iq(u)n, (i+1)q(u)n],~q(u)=u^{-1}v(u),$$ and
\BQNY
&&\Sigma\Sigma(u,n)\\
&&\coloneqq\sum_{i\neq j, 0\leq i,j\leq N(u,n)}\pk{\sup_{t\in I_i(n,u)}Q(t)>u, \sup_{t\in I_j(n,u)}Q(t)>u}\\
&&{=}\sum_{i\neq j, 0\leq i,j\leq N(u,n)}\pk{\sup_{t\in I_i(n,u), s\geq t}(B_{\alpha}(s)-B_{\alpha}(t)-c(s-t))>u, \sup_{t\in I_j(n,u), s\geq t}(B_{\alpha}(s)-B_{\alpha}(t)-c(s-t))>u}\\
&&=\sum_{i\neq j, 0\leq i,j\leq N(u,n)}\pk{\sup_{t\in I_i'(n,u), s\geq t}Z(s,t)>u^{1-\alpha/2}, \sup_{t\in I_j'(n,u), s\geq t}Z(s,t)>u^{1-\alpha/2}},
\EQNY
where in the last equality we use the self-similarity of fBm.
 Moreover, let
 $$ L_i(u)=[\tau^*+iq(u)n_1, \tau^*+(i+1) q(u) n_1], \quad M(u)=\left[\frac{u^{\alpha/2}\ln u}{v(u)n_1}\right],$$ $$ G(u)=\{s:|s-\tau^*|< u^{\alpha/2-1}\ln u\}, \quad G^c(u)=[0,\IF)\setminus G(u),$$
 and
 \BQN\label{pi22}\pi_2(u)=\sum_{-M(u)-1\leq i<j\leq M(u)+1} \pk{\sup_{t\in [0,q(u)n], s\in L_i(u)}Z(s,t)>u^{1-\alpha/2}, \sup_{t\in [0,q(u)n], s\in L_j(u)}Z(s,t)>u^{1-\alpha/2}}.\EQN

{\it \underline{Conditions A1 and A3.}}
Condition A1   follows from  Theorems 3.1-3.3 of \cite{KrzysPeng2015} while
A3 is due to Lemma 5.6 of \cite{KrzysPeng2015} and the upper bounds of $\Sigma_i(u), i=1,2,3,4$
in the proof of Theorem 3.1 in \cite{KrzysPeng2015}.

{\it \underline{Condition A2.}}
Due to stationarity of the process
$Q$, in order to show (\ref{Pickands})  it suffices
to find the exact asymptotics of $\pk{\int_{[0,v(u)n]}\I(Q(t)>u)dt>\vv(u)x}$ as $u\to\infty$.
By the self-similarity of $B_{\alpha}$, we have
\BQNY
\pk{\int_{[0,v(u)n]}\I(Q(t)>u)dt>\vv(u)x}&=&\pk{\int_{[0,v(u)n]}\I\left(\sup_{s\geq t}( B_{\alpha}(s)-B_{\alpha}(t)-c(s-t))>u\right)dt>\vv(u)x}\\
&=&\pk{\int_{[0, q(u)n]}\I\left(\sup_{s\geq t}Z(s,t)>u^{1-\alpha/2}\right)dt> q(u) x}.
\EQNY
\BEL\label{l.q}
For $n>x$
\BQN\label{shortQ}
\pk{\int_{[0, q(u)n]}\I\left(\sup_{s\geq t}Z(s,t)>u^{1-\alpha/2}\right)dt> q(u) x}\sim \widehat{\mathcal{B}}_{\alpha,\alpha}(x, n)\sqrt{\frac{2A}{B}}\frac{u}{m(u)v(u)}\Psi(u), \quad u\rw\IF.
\EQN
\EEL
{\it Proof.}
 {\it\underline{Upper bound}}.
Using the fact that
 $$\I\left(\sup_{s\geq t}Z(s,t)>u^{1-\alpha/2}\right)\leq \I\left(\sup_{s\in G(u) }Z(s,t)>u^{1-\alpha/2}\right)+\I\left(\sup_{s\in G^c(u)}Z(s,t)>u^{1-\alpha/2}\right)$$
we obtain
\BQNY
&&\pk{\int_{[0, q(u)n]}\I\left(\sup_{s\geq t}Z(s,t)>u^{1-\alpha/2}\right)dt> q(u) x}\\
&&\quad\leq \pk{\int_{[0, q(u) n]}\left(\I\left(\sup_{s\in G(u) }Z(s,t)>u^{1-\alpha/2}\right)+\I\left(\sup_{s\in G^c(u) }Z(s,t)>u^{1-\alpha/2}\right)\right)dt> q(u) x}\\
&&\quad\leq \pk{\int_{[0,q(u)n]}\I\left(\sup_{s\in G(u) }Z(s,t)>u^{1-\alpha/2}\right)dt> q(u) x}\\
&&\quad\quad +\pk{\int_{[0,q(u)n]}\I\left(\sup_{s\in G^c(u) }Z(s,t)>u^{1-\alpha/2}\right)dt> q(u) x}\\
&&\quad \quad+ \pk{\int_{[0,q(u)n]}\I\left(\sup_{s\in G(u) }Z(s,t)>u^{1-\alpha/2}\right)dt>0, \int_{[0, q(u)n]}\I\left(\sup_{s\in G^c(u) }Z(s,t)>u^{1-\alpha/2}\right)dt>0}\\
&&\quad \leq \pk{\int_{[0,q(u)n]}\I\left(\sup_{s\in G(u) }Z(s,t)>u^{1-\alpha/2}\right)dt> q(u) x}+2\pk{\sup_{t\in [0,q(u) n], s\in G^c(u)}Z(s,t)>u^{1-\alpha/2}}.
\EQNY
Moreover, since
\BQNY\label{upper}\I\left(\sup_{s\in G(u) }Z(s,t)>u^{1-\alpha/2}\right)\leq \sum_{|i|\leq M(u)+1} \I\left(\sup_{s\in L_i(u) }Z(s,t)>u^{1-\alpha/2}\right)
\EQNY
we have
\BQNY
\pk{\int_{[0,q(u) n]}\I\left(\sup_{s\in G(u) }Z(s,t)>u^{1-\alpha/2}\right)dt> q(u) x}\leq\pi_1(u)+\pi_2(u),
\EQNY
 where $\pi_2(u)$ is given in (\ref{pi22}) and
 \BQN\label{pi2}
\pi_{1}(u)=\sum_{|i|\leq M(u)+1}\pk{\int_{t\in [0,q(u)n]}\I\left(\sup_{s\in L_i(u) }Z(s,t)>u^{1-\alpha/2}\right)dt> q(u) x}.
 \EQN

By  Lemma 5.6 of \cite{KrzysPeng2015} we obtain  
 $$\pk{\sup_{t\in [0,q(u)n], s\in G^c(u)}Z(s,t)>u^{1-\alpha/2}}=o\left( \pk{\sup_{t\in [0,v(u)n]}Q(t)>u}\right), \quad u\rw\IF,$$
 and
in light of  the upper bounds of $\Lambda_i(u), i=1,2,3,4$ in the proof of Theorem 3.1 of \cite{KrzysPeng2015}
\BQN\label{pii2}
\pi_2(u)= o\left(\pk{\sup_{t\in [0,v(u)n]}Q(t)>u}\right), \quad u\rw\IF, n_1\rw\IF.
\EQN
Next we focus on $\pi_1(u)$.  We denote
$$ m(u)=\frac{1+c\tau^*}{(\tau^*)^{\alpha/2}}u^{1-\alpha/2},\quad
\tau^*=\frac{\alpha}{c(2-\alpha)},  $$
$$ A=\left(\frac{\alpha}{c(2-\alpha)}\right)^{-\alpha/2}\frac{2}{2-\alpha}, \quad  B=\left(\frac{\alpha}{c(2-\alpha)}\right)^{-\alpha/2-1}\frac{\alpha}{2}.$$
Rewrite
\BQNY
\pk{\int_{t\in [0,q(u)n]}\I\left(\sup_{s\in L_i(u) }Z(s,t)>u^{1-\alpha/2}\right)dt> q(u) x}=\pk{\int_{t\in [0,n]}\I\left(\sup_{s\in [0,n_1]}Z_{u,i}(s,t)>m(u) \right)dt>x},
\EQNY
where
$$Z_{u,i}(s,t)=\frac{B_{\alpha}(\tau^*+ q(u) (in_1+s))-B_{\alpha}( q(u) t)}{1+c(\tau^*+ q(u) (in_1+s-t))}\cdot\frac{1+c\tau^*}{(\tau^*)^{\alpha/2}}.$$
Let for $0<\epsilon<1$
$$m_i^{\pm}(u)=m(u)\left(1+\left(\frac{B}{2A}\pm \epsilon\right) q(u) (in_1\pm n )^2\right).$$
A direct calculation shows (see also Lemmas 5.3-5.4 in \cite{KrzysPeng2015}) that
\BQN\label{varaince}
m_i^-(u)\leq m(u)(Var(Z_{u,i}(s,t)))^{-1/2}\leq m_i^+(u), \quad |i|\leq M(u)+1
\EQN
and
\BQN\label{correlation}
\lim_{u\rw\IF}\sup_{|i|\leq M(u)+1}\sup_{(s,t)\neq (s',t'), (s,t), (s',t')\in [0,n_1]\times[0,n]}\left|(m_i^{\pm}(u))^2\frac{1-Corr(Z_{u,i}(s,t), Z_{u,i}(s',t'))}{|t-t'|^{\alpha}+|s-s'|^{\alpha}}-1\right|=0.
\EQN
Hence
\BQNY
\pk{\int_{t\in [0,n]}\I\left(\sup_{s\in [0,n_1]}Z_{u,i}(s,t)>m(u)\right)dt>x}\leq \pk{\int_{t\in [0,n]}\I\left(\sup_{s\in [0,n_1]}\overline{Z}_{u,i}(s,t)>m_i^{-}(u)\right)dt>x}.
\EQNY

Next, by
Lemma \ref{the-weak-conv}
applied to
$\Gamma: C([0,n]\times[0,n_1])\rw C([0,n])$ defined by
$\Gamma(f)=\sup_{t\in [0,n_1]}f(s,t),  f\in C([0,n]\times[0,n_1])$,
with $h=0$ in {\bf C0-C1}
and
{\bf C2} satisfied with $\zeta(s,t)=B_{\alpha}(s)+B_{\alpha}'(t)$,
we have
\BQN\label{uniformqueue}
\lim_{u\rw\IF}\sup_{|i|\leq M(u)+1}\left|\frac{\pk{\int_{t\in [0,n]}\I\left(\sup_{s\in [0,n_1]}\overline{Z}_{u,i}(s,t)>m_i^{-}(u)\right)dt>x}}{\Psi(m_i^-(u))}-\widehat{\mathcal{B}}_{\alpha,\alpha}(x, [0,n]\times[0,n_1])\right|=0,
\EQN
 and in light of Lemma \ref{lemA0}, we have
\BQN\label{pi1}
\pi_{{1}}(u)&\leq& \widehat{\mathcal{B}}_{\alpha,\alpha}(x, [0,n]\times[0,n_1])\sum_{|i|\leq M(u)+1}\Psi(m_i^-(u))\nonumber\\
 &\leq&\widehat{\mathcal{B}}_{\alpha,\alpha}(x, [0,n]\times[0,n_1])\Psi(u)\sum_{|i|\leq M(u)+1}e^{-m^2(u)\left(\frac{\rrd{B}}{2A}- \epsilon\right)\left(u^{-1}v(u)(in_1)\right)^2}\nonumber\\
 &\leq&\frac{\widehat{\mathcal{B}}_{\alpha,\alpha}(x, [0,n]\times[0,n_1])}{n_1}\sqrt{\frac{2A\pi}{B}}\frac{u}{m(u)v(u)}\Psi(u)\nonumber\\
 &\sim &\widehat{\mathcal{B}}_{\alpha,\alpha}(x, n)\sqrt{\frac{2A\pi}{B}}\frac{u}{m(u)v(u)}\Psi(u),
\EQN
as $u\rw\IF, n_1\rw\IF, \epsilon\rw 0$. Therefore, we conclude that
\BQNY\label{upper4}
\pk{\int_{[0, q(u)n]}\I\left(\sup_{s\geq t}Z(s,t)>u^{1-\alpha/2}\right)dt> q(u) x}\leq \widehat{\mathcal{B}}_{\alpha,\alpha}(x, n)\sqrt{\frac{2A\pi}{B}}\frac{u}{m(u)v(u)}\Psi(u), \quad u\rw\IF.
\EQNY

{\it \underline{Lower bound}}. Observe that for u sufficiently large, $s>t$ holds for all $s\in G(u), t\in [0,q(u)n]$. Therefore,
\BQNY
\pk{\int_{[0, q(u)n]}\I\left(\sup_{s\geq t}Z(s,t)>u^{1-\alpha/2}\right)dt> q(u) x}\geq \pk{\int_{[0,q(u)n]}\I\left(\sup_{s\in G(u) }Z(s,t)>  u^{1-\alpha/2}\right)dt>
	q(u) x}.
\EQNY
By the fact that
\BQNY\label{lower}\I\left(\sup_{s\in G(u) }Z(s,t)>u^{1-\alpha/2}\right)&\geq& \sum_{|i|\leq M(u)} \I\left(\sup_{s\in L_i(u) }Z(s,t)>u^{1-\alpha/2}\right)\nonumber\\
&&-\sum_{-M(u)\leq i<j\leq M(u)} \I\left(\sup_{s\in L_i(u) }Z(s,t)>u^{1-\alpha/2}, \sup_{s\in L_j(u) }Z(s,t)>u^{1-\alpha/2}\right)\\
&=:& A_1(u,t)-A_2(u,t),\nonumber
\EQNY
it follows that for $\epsilon>0$ (recall $q(u)= u^{-1} v(u)$)
\BQN\label{lower1}
&&\pk{\int_{[0, q(u)n]}\I\left(\sup_{s\geq t}Z(s,t)>u^{1-\alpha/2}\right)dt> q(u) x}\nonumber\\
&&\geq
\pk{\int_{[0,q(u)n]}\left(A_1(u,t)-A_2(u,t)\right)dt>  q(u) x}\nonumber\\
&&\geq \pk{\int_{[0,q(u)n]}A_1(u,t)dt> q(u)(x+\epsilon), \int_{[0,q(u)n]}A_2(u,t)dt<q(u)\epsilon}\nonumber\\
&&\geq \pk{\int_{[0, q(u) n]}A_1(u,t)dt>q(u)(x+\epsilon)}-\pk{\int_{[0, q(u)n]}A_2(u,t)dt\geq  q(u)  \epsilon}\nonumber\\
&&\geq \pk{\exists \rrd{i:}\, |i|\leq M(u), \int_{[0,q(u)n]} \I\left(\sup_{s\in L_i(u) }Z(s,t)>u^{1-\alpha/2}\right)dt> q(u)(x+\epsilon)}-\pi_2(u)\nonumber\\
&&\geq \sum_{|i|\leq M(u)}\pk{\int_{t\in [0, q(u) n]}\I\left(\sup_{s\in L_i(u) }Z(s,t)>u^{1-\alpha/2}\right)dt>q(u)(x+\epsilon)}-2\pi_2(u),
\EQN
where $\pi_2(u)$ is defined in (\ref{pi22}).
Similarly as in (\ref{pi1}) and in light of (\ref{pii2}), we have
\BQNY
\pk{\int_{[0, q(u)n]}\I\left(\sup_{s\geq t}Z(s,t)>u^{1-\alpha/2}\right)dt> q(u) x}
 &\geq& \widehat{\mathcal{B}}_{\alpha,\alpha}(x+\epsilon, n)\sqrt{\frac{2A}{B}}\frac{u}{m(u)v(u)}\Psi(u)\\
 &\geq& \widehat{\mathcal{B}}_{\alpha,\alpha}(x, n)\sqrt{\frac{2A}{B}}\frac{u}{m(u)v(u)}\Psi(u), \quad u\rw\IF,\epsilon\rw 0.
\EQNY
Consequently for $n>x$
\BQN\label{shortQ}
\pk{\int_{[0, q(u)n]}\I\left(\sup_{s\geq t}Z(s,t)>u^{1-\alpha/2}\right)dt> q(u) x}\sim \widehat{\mathcal{B}}_{\alpha,\alpha}(x, n)\sqrt{\frac{2A}{B}}\frac{u}{m(u)v(u)}\Psi(u), \quad u\rw\IF.
\EQN
\QED

Moreover, by Lemma \ref{lemA0}
\BQNY
\frac{\mathcal{B}_{\alpha}(x)}{\mathcal{B}_{\alpha}(0)}=\frac{\widehat{\mathcal{B}}_{\alpha,\alpha}(x)}{\widehat{\mathcal{B}}_{\alpha,\alpha}(0)}=\lim_{n\rw\IF}\frac{\widehat{\mathcal{B}}_{\alpha,\alpha}(x,n)}{\widehat{\mathcal{B}}_{\alpha,\alpha}(0,n)}\in (0,{1]}.
\EQNY
 Thus A2 holds with
$$ \bar F(x)=\frac{\mathcal{B}_{\alpha}(x)}{\mathcal{B}_{\alpha}(0)}, \quad x\geq 0.$$
This completes the proof of case ii).\\

For case i), note that if $x=0$, the claim clearly holds. Next we suppose that $0<x<T$. By (\ref{shortQ}) for any $0<\epsilon<\min(x/2, (T-x)/2)$,
\BQNY \pk{\int_{[0,T_u]}\I(Q(t)>u)dt>\vv(u)x}&\leq& \pk{\int_{[0,v(u)(T+\epsilon)]}\I(Q(t)>u)dt>\vv(u)x}\\
&\leq &\pk{\int_{[0,v(u)T]}\I(Q(t)>u)dt>\vv(u)(x-\epsilon)}\\
&\sim& \frac{\widehat{\mathcal{B}}_{\alpha,\alpha}(x-\epsilon,T)}{\widehat{\mathcal{B}}_{\alpha,\alpha}(0,T)}\pk{\sup_{t\in [0,v(u)T]}Q(t)>u}, \quad u\rw\IF.
\EQNY
Analogously,
\BQNY \pk{\int_{[0,T_u]}\I(Q(t)>u)dt>\vv(u)x}
\geq \frac{\widehat{\mathcal{B}}_{\alpha,\alpha}(x+\epsilon,T)}{\widehat{\mathcal{B}}_{\alpha,\alpha}(0,T)}\pk{\sup_{t\in [0,v(u)T]}Q(t)>u}, \quad u\rw\IF.
\EQNY
In light of Remark \ref{lemma} we establish the claim by letting $\epsilon\rw 0$ in the above inequalities.
This completes the proof. \QED

\section{Appendix}\label{s.supl}


\prooflem{the-weak-conv}
For notational simplicity denote by $\rho_{u,j}$ the correlation function of the random field $\xi_{u,j}$.
Further set
$$ \chi_{u,j}(s) \coloneqq g_{u,j}(\overline\xi_{u,j}(s)-\rho_{u,j}(s,\vk{0})\overline\xi_{u,j}(\vk{0})),\quad s\in E_1
$$

and
$$f_{u,j}(s,y)\coloneqq\rd{y} \rho_{u,j}(s,\vk{0})-g^2_{u,j}\LT( 1-\rho_{u,j}(s,\vk{0}) \RT)-g^2_{u,j}\frac{1-\sigma_{u,j}(s)}{\sigma_{u,j}(s)}  ,\ s\in E_1 , y\in\R.$$
Conditioning on $\xi_{u,j}(\vk{0})$, by {\bf F1} and using that $\overline{\xi}_{u,j}(\vk{0})$ and $\overline{\xi}_{u,j}(s)-\ruk(s,\vk{0})\overline{\xi}_{u,j}(\vk{0})$ are mutually independent we obtain
\BQNY
&&\pk{ \int_{ E_{2} } \I\left\{\Gamma\left({g_{u,j}(\xi_{u,j}(s)-g_{u,j})}\right)(t)>0\right\} \eta(dt) >x}\\
&&=\frac{e^{-g^2_{u,j}/2}} { \sqrt{2\pi}g_{u,j}}
\int_{\R} \exp\LT(-y-\frac{y^2}{2g^2_{u,j}}\RT) \pk{ \int_{ E_{2} } \I\left\{\Gamma\left(g_{u,j}(\xi_{u,j}(s)-g_{u,j})\right)(t)>0\right\} \medt > x | \xi_{u,j}(\vk{0})=g_{u,j}+yg_{u,j}^{-1}} dy\\
&&= \frac{e^{-g^2_{u,j}/2}} { \sqrt{2\pi}g_{u,j}}
\int_{\R} \exp\LT(-y-\frac{y^2}{2g^2_{u,j}}\RT) \pk{ \int_{ E_{2}} \I\left\{\Gamma\left(\sigma_{u,j}(s)\left(\chi_{u,j}(s) + f_{u,j}(s,y)\right)\right)(t)>0 \right\} \eta(dt)
	>x } dy\\
&&=\frac{e^{-g^2_{u,j}/2}} { \sqrt{2\pi}g_{u,j}} \int_{\R} \exp\LT(-y-\frac{y^2}{2g^2_{u,j}}\RT) \mathcal{I}_{u,j}(y;x) dy,
\EQNY
where
$$\mathcal{I}_{u,j}(y;x)\coloneqq\pk{ \int_{ E_{2}} \I\left\{\Gamma\left(\sigma_{u,j}(s)\left(\chi_{u,j}(s) + f_{u,j}(s,y)\right)\right)(t)>0 \right\} \eta(dt) >x }.$$
 Noting that
$$\lim_{u\rw\IF}\sup_{j\in S_u}\left|\frac{\frac{e^{-g^2_{u,j}/2}} { \sqrt{2\pi}g_{u,j}}}{\Psi(g_{u,j})}-1\right|=0
$$
in order to show the claim it suffices to prove that
\BQN\label{aux-uni-conv}
\lim_{u\to\IF}\sup_{j\in S_u}\left| \int_{\R} \exp\LT(-y-\frac{y^2}{2g^2_{u,j}}\RT)  \mathcal{I}_{u,j}(y;x) dy - \MB^{\Gamma, h,\eta}_{\zeta}( x, E_2)\right| =0
\EQN
for all $x\geq 0$. In view of  {\bf C3} it follows that  that  for $u>u_0$
$$Var(\chi_{u,j}(\vk{s})-\chi_{u,j}(\vk{s}')) \leq g^2_{u,j} \E{\overline\xi_{u,j}(\vk{s})-\overline\xi_{u,j}(\vk{s'})}^2 \leq Q_1 \norm{\vk{s}-\vk{s'}}^\nu, \quad s, s'\in E_1 ,
$$
with $\nu>0$.  Further, by {\bf C0-C2} for each $y\in\mathbb{R}$
\BQN\label{uniformtrend} \limit{u} \sup_{j\in S_u, s\in E_1}\abs{f_{u,j}(s,y)-y +\sigma^2_{\zeta}(s)+h(s)}
=0.
\EQN
Hence, by {\bf F2}   
\BQN\nonumber
\sup_{j\in S_u}e^{-y}\iukxz &\leq& e^{-y}\sup_{j\in S_u}\pk{\sup_{t\in E_2}\Gamma\left(\chi_{u,j}(s) + f_{u,j}(s,y)\right)(t)>0} \\ \nonumber
&\leq&
e^{-y}\sup_{j\in S_u}\pk{\sup_{s\in E_1}\{\chi_{u,j}(s) + f_{u,j}(s,y)\}>0} \\ \nonumber
&\leq& e^{-y}\sup_{j\in S_u}\pk{\sup_{s\in E_1} \chi_{u,j}(s)>Q_2\abs{y}-Q_3} \\ \label{for-domi-Iuk}
&\leq& Q_4 \abs{y}^{2n/\nu-1}e^{-Q_5y^2-y},\quad y<-M,
\label{for-domi-IFa}
\EQN
where in the  last inequality we used Piterbarg inequality {and $M>0$}. Moreover,  it follows trivially that  for all $x\geq 0$ 
\BQN\label{upper1} \sup_{j\in S_u}e^{-y}\iukxz \leq e^{-y},\quad y\in\mathbb{R}.
\EQN
Therefore by the dominated convergence theorem and assumption {\bf C0}
\BQNY
&&\sup_{j\in S_u}\left|\int_{\R}  \exp\LT(-y-\frac{y^2}{2g^2_{u,j}} \RT)  \iukxz dy -\int_{\R}e^{-y}\iukxz dy \right|\\
&&\leq \int_{\R} \sup_{j\in S_u} \LT(e^{-y} \iukxz({1- e^{-y^2/(2g^2_{u,j})}  })\RT) dy \rightarrow 0, \quad u\to\IF.
\EQNY
Hence in order to prove the convergence in \eqref{aux-uni-conv} it suffices to show that
\BQN\label{aux-2-uni-conv}
\lim_{u\to\IF}\sup_{j\in S_u}\abs{\int_{\R} e^{-y} \mathcal{I}_{u,j}(y;x) dy - \MB^{\Gamma,h,\eta}_{\zeta}(x, E_2) } =0
\EQN
for all $x\in [0,\eta(E_2))$.\\
{\it \underline{Weak convergence}}. The claim follows from the same arguments as in \cite{MR4127347}[Lem 4.3,4.7],  where the precise meaning of uniform weak convergence is also given. Thus
let  $C( E_1 )$ denote the Banach space of all continuous functions on the compact set $ E_1 $ equipped with supremum norm.
For any $\vk{s},\vk{s}'\in E_1 $, by {\bf C2} we have
$$Var(\chi_{u,j}(s)-\chi_{u,j}(\vk{s}')) = g^2_{u,j}\LT( \E{\overline\xi_{u,j}(s)-\overline\xi_{u,j}(\vk{s}')}^2 - \LT( \rho_{u,j}(s,\vk{0})-\rho_{u,j}(\vk{s}',\vk{0}) \RT)^2 \RT) \rightarrow 2Var(\zeta(s)-\zeta(\vk{s}'))  $$
uniformly with respect to $j\in  S_u$ as $u\to\IF$.
Hence, the finite-dimensional distributions of $\chi_{u,j}(s), s\in E_1$ weakly converge to that
of $\sqrt{2}\zeta(s),s\in E_1 $ uniformly with respect to $j\in  S_u$. In view of {\bf C3},
we know that the measures on $C( E_1 )$ induced by $\{\chi_{u,j}(s),s\in E_1 ,j\in S_u\}$ are uniformly tight for large $u$, and by {\bf C1}, $\sigma_{u,j}(s)$ converges to $1$ uniformly for $s\in  E_1 $ and $j\in  S_u$ as $u\to \IF$.
Therefore,  $\{\sigma_{u,j}(s)\chi_{u,j}(s),s\in E_1 \}$
converge weakly to $\{\sqrt{2}\zeta(s),s\in E_1 \}$ as $u\to\IF$ uniformly with respect to $j\in S_u$, which together with (\ref{uniformtrend})
implies  that for each $y\in \mathbb{R}$, the probability measures on $C( E_1 )$ induced by $\{\chi_{u,j}^f(s,y),s\in E_1\}$ converges weakly as $u\to\IF$ to that induced by $\{\zeta_h(s) + y,\vk{t}\in E_1 \}$ uniformly with respect to $j\in S_u$, where
$$\chi_{u,j}^f(\vk{s},y)\coloneqq\sigma_{u,j}(s)\left(\chi_{u,j}(s) + f_{u,j}(s,y)\right)\quad \textrm{and}\quad \zeta_h(s)\coloneqq\sqrt{2}\zeta(s)-\sigma^2_{\zeta}(t)-h(s).$$
Continuous mapping theorem implies that
for each $y\in \mathbb{R}$, the \cEE{push-forward} probability measures \cEE{$P_{u,y}$} on $C( E_{2} )$ induced by
$\{\Gamma\left(\chi_{u,j}^f(\cdot,y)\right)(t),t\in E_{2}\}$ converges weakly  the push-forward probability measure \cEE{$P_{y}$}  induced by $\{\Gamma\left(\zeta_h\right)(t) + y,t\in E_{2} \}$ as $u\to\IF$ uniformly with respect to $j\in S_u$.

The continuity of the sojourn functional is also discussed in \cite{berman1973excursions}[Lem 4.2].
A sequence of functions $f_n \in C(E_2)$ converges to $f\in C(E_2)$ as $n \to \IF$ with respect to uniform topology if $f_n \to f$ uniformly as $n\to \IF$.  Since $\eta$ is absolutely continuous with respect to Lebesgue measure on $E_2$ we can define the set $$A_*=\left\{ f\in C( E_{2} ): \int_{ E_{2} } \mathbb{I}(f(t)=0)\eta(dt)>0\right\},$$
\cEE{which is  measurable in the completion $\mathcal{C}^\mu$ of  $\mathcal{C}$ with respect to $\nu$, where $\mathcal{C}$ is} the Borel $\sigma$-field of $C_2(E)$. Its complement  belongs to $\cEE{\mathcal{C}^\mu}$, i.e.,
$$A_*^c= C(E_2) \setminus A_* \in \cEE{\mathcal{C}^\mu}.$$
 Any function $f\in A_*^c$ is a continuity point of
the sojourn functional $J: C(E_2) \mapsto [0, \eta(E_2)],$ where
$$J(f)= \int_{E_2} 1(f(t)> 0)\eta(dt), f\in C(E_2).$$
This functional is measurable $\mathcal{C}/\mathcal{B}(\R)$ by the assumption on $\eta$.
We shall show that it is continuous at any $f\in A_*^c$.
Let such $f$ be given. By the definition of the integral such $f$ is not equal to zero on any compact interval of $\R$.
Let $f_n \to f$ uniformly as $n\to \IF$. Then $1(f_n(t)> 0) \to 1(f(t)> 0)$ as $n\to \IF$
for almost all $t\in \R$ (with respect to Lebesgue measure).
Hence by dominated convergence theorem we have $J(f_n) \to J(f)$ as $n\to \IF$, which means that the functional is continuous for all $f\in A_*^c$. Recall that $P_y$ \cEE{is the push-forward (image measure)} on $C(E_2)$ with respect to $\Gamma( \xi_h)+ y$.
We claim that
$$ P_y(A_*)>0$$
is possible only for $y$ in a countable set of $\R$.
Indeed, any $f\in A_*$ is such that it is constant equal to zero on a compact interval.
Consequently, $P_y(A_*)>0$ means that the functions $f\in A_*$ are constant equal to $ -y$ on some interval of $\R$.
If this is true for two different $y$'s, then the intervals where $f$ is constant equal $-y$ must be disjoint,
therefore this can be true only for countable $y$'s.

Alternatively,
using the fact that $\pk{\Gamma(\zeta_h)(t)+y=0}=0$ a.e., $y\in\mathbb{R}$,
by  the $\sigma$-finiteness of $\eta$,  Fubini-Tonelli theorem yields
\BQNY
\int_\mathbb{R}\E{\int_{ E_{2}} \mathbb{I}(\Gamma(\zeta_h)(t)+y=0)\eta(dt)}dy=\int_{ E_{2}} \int_\mathbb{R}\pk{\Gamma(\zeta_h)(t)+y=0}dy\eta(dt)
=0.
\EQNY
Hence  for almost all $y\in \mathbb{R}$
$$\E{\int_{ E_{2}} \mathbb{I}(\Gamma(\zeta_h)(t)+y=0)\eta(dt)}=0,$$
which means that,   for almost all $y\in \mathbb{R}$
$$P_y(A_*)=\mathbb{P}\left(\int_{ E_{2}} \mathbb{I}(\Gamma(\zeta_h)(t)+y=0)\eta(dt)>0\right)=0.
$$
Consequently, since $J(f)$ is continuous for $f\in A_*^c$,  by continuous mapping theorem, as $u\to\infty$
\BQN\label{wea-conv-cont-aux}\int_{ E_{2} } \mathbb{I}\left(\Gamma\left(\chi_{u,j}^f(\cdot,y)\right)(t)>0\right)\eta(dt)
\EQN
weakly converges  to
$$\int_{ E_{2} } \mathbb{I}\left(\Gamma\left(\zeta_h\right)(t)+y>0\right)\eta(dt)$$
uniformly with respect to $j\in S_u$ for almost all $y\in \mathbb{R}$.\\
{\it \underline{Convergence on continuity points.}}
Define
$$\mathcal{I}(y;x) \coloneqq \pk{ \int_{ E_{2} } \mathbb{I}\left(\Gamma(\zeta_h)(t) + y >0\right) \medtl >x}. $$
We draw a similar argument as in Theorem 1.3.1 of \cite{Berman92} to verify \eqref{aux-2-uni-conv} for all  continuity points $x\in (0,\eta(E_2))$ of  $\MB^{\Gamma, h,\eta}_{\zeta}(x, E_2)$.
Let $x_0\in (0, \eta(E_2))$ be such a continuity point, that is
$$ \lim_{\vp\to0}\int_{\R} \LT(\mathcal{I}(y;x_0+\vp) -\mathcal{I}(y;x_0-\vp)\RT)e^{-y} dy=0.$$
Since for large $M$ and all $x\geq 0$ by {\bf F2}  as in the derivation of \eqref{for-domi-IFa}  we have
\BQN\label{for-domi-IF}
e^{-y}\mathcal{I}(y;x) \leq Q_4' \abs{y}^{2n/\nu-1}e^{-Q_5y^2-y},\quad y<-M
\EQN
it follows from the dominated convergence theorem that
$$ \int_{\R} \LT(\mathcal{I}(y;x_0+) -\mathcal{I}(y;x_0-)\RT) e^{-y} dy=0
$$
and thus by the monotonicity of $\mathcal{I}(y;x)$ in $x$ for each fixed $y$,
$x_0$ is a continuous point of $\mathcal{I}(y;x)$ for a.e. $y\in\R$.
Thus by \eqref{wea-conv-cont-aux} for a.e. $y\in\R$
\BQN\label{exp-almost-z}
\lim_{u\to\IF}\sup_{j\in  S_u}\abs{ \mathcal{I}_{u,j}(y;x_0)  - \mathcal{I}(y;x_0)} =0.
\EQN
As shown in \eqref{for-domi-Iuk}, (\ref{upper1}) and \eqref{for-domi-IF}
it follows from the dominated convergence theorem that
\BQN\label{new1}
&&\sup_{j\in S_u}\abs{\int_{\R} e^{-y} \mathcal{I}_{u,j}(y;x_0) dy - \int_{\R}e^{-y} \mathcal{I}(y;x_0)dy } \nonumber\\
&&\leq\int_{\R} \sup_{j\in S_u}\abs{\mathcal{I}_{u,j}(y;x_0) -  \mathcal{I}(y;x_0)} e^{-y} dy  \rightarrow 0,\quad u\to\IF
\EQN
establishing the proof for all continuity points $x\in (0,\eta(E_2))$. Moreover, for the case that $x=0$, (\ref{new1}) also holds by replacing sojourn with supremum. This can be shown directly without any continuity requirement for $\MB^{\Gamma, h,\eta}_{\zeta}(x, E_2)$ at $x=0$.\\
{\it \underline{Continuity of $\MB^{\Gamma, h,\eta}_{\zeta}(x, E_2)$}}.
Next we show that $\MB^{\Gamma, h,\eta}_{\zeta}(x, E_2)$ is  continuous at any  $x\in (0, \eta(E_2))$ using that $\eta$ is equivalent with
Lebesgue measure on $E_2$.  Note that  $\MB^{\Gamma, h,\eta}_{\zeta}(x, E_2)$ is clearly right continuous at $0$. Next we show the continuity at $x \in (0, E_{2})$. The claimed continuity at $x$ follows if we show
$$\int_\R   \pk{A_y} 
e^{-y} dy=0, \quad A_y=\left\{ \int_{ E_{2} } \mathbb{I}\big(\Gamma(\zeta_h)(t) + y >0\big) \medtl =x \right\}, \quad y\in\mathbb{R}.$$
If  $$\int_{ E_{2} } \mathbb{I}\big(\Gamma(\zeta_h)(t) + y >0\big) \medtl =x,$$
with $0<x<\eta(E_{2})$, then using the fact that $\Gamma(\zeta_h)(t)$ is continuous over $ E_{2}$ and the Lebesgue measure is
absolutely continuous with respect to $\eta$,  we have that for any $y'> y$\\
$$\int_{ E_{2} } \mathbb{I}\big(\Gamma(\zeta_h)(t) + y' >0\big) \medtl >x.$$
This implies that
$
A_y\cap A_{y'}=\emptyset, y\neq y', y,y'\in \mathbb{R}.
$
Noting that the continuity of $\Gamma(\zeta_h)$ guarantees the measurability of $A_y$, and
$$\{y: y\in \mathbb{R} \quad \text{such that} \quad \pk{A_y}>0\}$$
is a countable set because if it were not we would find countably many (disjoint) $A_y$ such that $\sum\pk{A_y}=\infty$.

Thus we get $\int_\R  \pk{A_y} e^{-y} dy=0, $ hence $\MB^{\Gamma,h,\eta}_{\zeta}(x, E_2)$ is continuous on $(0,\eta(E_{2}))$, establishing the claim. 
\QED

Before proceeding to the proof of Lemma \ref{l.1},
under notation introduced in the proof of Proposition \ref{THSTA},
we denote and analyze
\BQN\Sigma\Sigma_1(u,n)&\coloneqq&\sum_{0\leq k_i, k_i'\leq N_i(u,n), i=1,2, (k_1,k_2)\neq (k_1',k_2')}\pk{\sup_{t\in I_{k_1,k_2}(u,n)}X(t)>u, \sup_{t\in I_{k_1',k_2'}(u,n)}X(t)>u },\label{ssigma1}\\
\Sigma\Sigma_2(u,n)&\coloneqq&\sum_{0\leq 2k_i, 2k_i'\leq N_i(u,n), i=1,2, (k_1,k_2)\neq (k_1',k_2')}\pk{\sup_{t\in I_{2k_1,2k_2}(u,n)}X(t)>u, \sup_{t\in I_{2k'_1,2k'_2}(u,n)}X(t)>u },\label{ssigma2}\\
 \Theta(u)&\coloneqq&T_1T_2a_1^{1/\alpha_1}a_2^{1/\alpha_2}u^{2/\alpha_1+2/\alpha_2}\Psi(u).
\label{theta}
\EQN

Moreover, following notation introduced in the proof of Proposition \ref{Onepoint},
let
$$\Sigma\Sigma_3''(u,n)\coloneqq \sum_{|k_i|, |k_i'|\leq N_i'(u,n), i=1,2, (k_1,k_2)\neq (k_1',k_2')}\pk{\sup_{t\in I_{k_1,k_2}(u,n)}X(t)>u, \sup_{t\in I_{k_1',k_2'}(u,n)}X(t)>u }.$$
\BQN\label{E1}
 \hat{I}_{k_2}(u,n)\coloneqq I_{-1,k_2}(u,n)\cup I_{0,k_2}(u,n),\quad  E_1(u,n)\coloneqq\bigcup_{|k_2|\leq N_2'(u,n)}\hat{I}_{k_2}(u,n),
 \EQN
 and
 \BQN\Sigma_3'(u,n)&\coloneqq&\sum_{ |k_i|\leq N_i'(u,n)+1, i=1,2,~ k_1\neq -1, 0}\pk{\sup_{t\in I_{k_1,k_2}(u,n)}X(t)>u}, \label{eq3}\\
 \Sigma\Sigma_3(u,n)&\coloneqq&\sum_{|k_2|, |k_2'|\leq N_2'(u,n), k_2\neq k_2'}\pk{\sup_{t\in \hat{I}_{k_2}(u,n)}X(t)>u, \sup_{t\in \hat{I}_{k_2'}(u,n)}X(t)>u }, \label{eq4}\\
 \Sigma\Sigma_4(u,n)&\coloneqq&\sum_{ |2k_2|,|2k_2'|\leq N_2'(u,n)-1, k_2\neq k_2'}\pk{\sup_{t\in \hat{I}_{2k_2}(u,n)}X(t)>u, \sup_{t\in \hat{I}_{2k_2'}(u,n)}X(t)>u }.\label{ssigma4}
 \EQN

\BEL\label{gaussianfields}
Under the assumptions of Proposition \ref{THSTA}
\BQN\label{Stationary}\pk{ \sup_{t\in E} X(t)>u}\sim \sum_{0\leq k_i\leq N_i(u,n), i=1,2}\pk{ \sup_{t\in I_{k_1,k_2}(u,{n})} X(t)>u}\sim \CC_0\Theta(u), \quad u\rw\IF,
n\to\infty,
\EQN
where $\CC_0>0$.
Moreover,  for all large $u$ and $n$
\BQNY
\Sigma\Sigma_1(u,n)&\leq&\left(\frac{\CC_2}{\sqrt{n}}+ e^{-\CC_1n^{\CC}}\right)\Theta(u),\quad \Sigma\Sigma_2(u,n)\leq e^{-\CC_1n^{\CC}}\Theta(u),
\EQNY
where $\CC,\CC_1$ and $\CC_2$ are some positive constants.
\EEL
\prooflem{gaussianfields}
Asymptotics (\ref{Stationary})  follow from Lemma 7.1 in \cite{Pit96},
while the bounds can be deduced from equations (7.4) and (7.6) in the proof of Lemma 7.1 in \cite{Pit96}. \QED

\BEL\label{gaussianfields1} Under the assumptions
of Proposition \ref{Onepoint}, for
	  $\alpha_i<\beta_i, i=1,2$, 
$$\pk{\sup_{t\in E\setminus E(u,n)}X(t)>u}=o\left(\pk{\sup_{t\in E}X(t)>u}\right)$$
as $u\to\infty, n\to\infty$, and
$$\Sigma\Sigma_3''(u,n)=o\left(\sum_{0\leq k_i\leq N_i'(u,n), i=1,2}\pk{\sup_{t\in I_{k_1,k_2}(u,n)}X(t)>u}\right),$$
as $u\to\infty, n\to\infty$.
For $\alpha_1=\beta_1, \alpha_2<\beta_2$
$$\pk{\sup_{t\in E\setminus E_1(u,n)}X(t)>u}=o\left(\pk{\sup_{t\in E}X(t)>u}\right),$$
as $u\to\infty$, $n\to\infty$,
and for $u$ and $n$ sufficiently large
\begin{align*}
\Sigma\Sigma_3(u,n)&\leq \left(\frac{\CC_2}{\sqrt{n}}+e^{-\mathbb{C}_1n^{\mathbb{C}}}\right)\pk{\sup_{t\in E}X(t)>u},\\
\Sigma_3'(u,n)&\leq e^{-\mathbb{C}_1n^{\mathbb{C}}}\pk{\sup_{t\in E}X(t)>u},\\
\Sigma\Sigma_4(u,n)&\leq e^{-\mathbb{C}_1n^{\mathbb{C}}}\pk{\sup_{t\in E}X(t)>u}.
\end{align*}
For $\alpha_1=\beta_1$ and $\alpha_2=\beta_2$ 
$$\pk{\sup_{t\in E\setminus \bigcup_{i,j\in \{-1,0\}}I_{i,j}(u,n)}X(t)>u}=o\left(\pk{\sup_{t\in E}X(t)>u}\right),$$
as $u\to\infty, n\to\infty$.
\EEL
\prooflem{gaussianfields1}
The proof of Lemma \ref{gaussianfields1} follows from \cite{KEP20151}. Specifically,
the first one follows from  (34),  the second one from (40) and (41),
the third one from (34) and (46), the fourth one from (48) and (49),
the fifth one from (46),   the six one from (48),
and the last one from (34) and (52) in the proof of  Theorem 3.1 of \cite{KEP20151}.  \QED
Now we are in the position to prove Lemma \ref{l.1}.
\\

\prooflem{l.1}
{\underline{ \it Ad (i)}.} We follow notation introduced in the proof of Proposition \ref{THSTA}.
For any $n,n_1>\sqrt{x}$, we have
\BQN\label{bonferroni}
\Sigma_1^-(u,n_1)-\Sigma\Sigma_1(u,n_1)\leq \pk{\int_{E(u,n)}\mathbb{I}(X(t)>u)dt >\vv(u)x}\leq \Sigma_1^+(u,n)+\Sigma\Sigma_1(u,n),
\EQN
where $\Sigma\Sigma_1(u,n)$ is given in (\ref{ssigma1}) and
\BQNY
\Sigma_1^\pm(u,n)=\sum_{0\leq k_i\leq N_{i}(u,n)\pm 1, i=1,2 }\pk{\int_{I_{k_1,k_2}(u,n)}\mathbb{I}(X(t)>u)dt >\vv(u)x}.
\EQNY
By (\ref{uniform1}), it follows that
\BQNY
\Sigma_1^+(u,n)&\leq& \sum_{0\leq k_i\leq N_{i}(u,n), i=1,2 }\mathcal{B}_{\alpha_1,\alpha_2}(x, [0,n]^2)\Psi(u)\\
&\leq& \frac{\mathcal{B}_{\alpha_1,\alpha_2}(x, [0,n]^2)}{n^2} \Theta(u), \quad u\rw\IF,
\EQNY
where $\Theta(u)$ is defined in (\ref{theta}).
Analogously, we obtain the lower bound
\BQNY
\Sigma_1^-(u,n)
&\geq& \frac{\mathcal{B}_{\alpha_1,\alpha_2}(x, [0,n]^2)}{n^2} \Theta(u), \quad u\rw\IF.
\EQNY
Lemma \ref{gaussianfields} shows that for $u$ and $n$ sufficiently large
$$\Sigma\Sigma_1(u,n)\leq\left(\frac{\CC_2}{\sqrt{n}}+ e^{-\CC_1n^{\CC}}\right)\Theta(u).
$$
Dividing both sides of (\ref{bonferroni}) by $\Theta(u)$ and letting $u\rw\IF$, we have
\BQNY
\frac{\mathcal{B}_{\alpha_1,\alpha_2}(x, [0,n_1]^2)}{n_1^2}-\frac{\CC_2}{\sqrt{n_1}}-e^{-\mathbb{C}_1n_1^{\mathbb{C}}}\leq \frac{\mathcal{B}_{\alpha_1,\alpha_2}(x, [0,n]^2)}{n^2}+\frac{\CC_2}{\sqrt{n}}+ e^{-\mathbb{C}_1n^{\mathbb{C}}}.
\EQNY
The above implies that
\BQNY
\limsup_{n\rw\IF}\frac{\mathcal{B}_{\alpha_1,\alpha_2}(x, [0,n]^2)}{n^2}&=& \liminf_{n\rw\IF}\frac{\mathcal{B}_{\alpha_1,\alpha_2}(x, [0,n]^2)}{n^2}<\IF.
\EQNY
Next we show that
$$\liminf_{n\rw\IF}\frac{\mathcal{B}_{\alpha_1,\alpha_2}(x, [0,n]^2)}{n^2}>0.$$
Observe that
\BQN\label{lowerfields}
\pk{\int_{E}\mathbb{I}(X(t)>u)dt >\vv(u)x}\geq \Sigma_2(u,n)-\Sigma\Sigma_2(u,n),
\EQN
where $\Sigma\Sigma_2(u)$ is given in (\ref{ssigma2}) and
\BQNY
\Sigma_2(u,n)=\sum_{0\leq 2k_i\leq N_{i}'(u,n), i=1,2 }\pk{\int_{I_{2k_1,2k_2}(u,n)}\mathbb{I}(X(t)>u)dt >\vv(u)x}.
\EQNY
In light of (\ref{uniform1}), we have
\BQNY
\Sigma_2(u,n)&\geq& \sum_{0\leq 2k_i\leq N_{i}'(u,n), i=1,2 }\mathcal{B}_{\alpha_1,\alpha_2}(x, [0,n]^2)\Psi(u)\\
&\geq& \frac{\mathcal{B}_{\alpha_1,\alpha_2}(x, [0,n]^2)}{4n^2} \Theta(u), \quad u\rw\IF.
\EQNY
Moreover, by Lemma \ref{gaussianfields} we have, for $u$ and $n$ large enough
$$\Sigma\Sigma_2(u,n)\leq e^{-\CC_1n^{\CC}}\Theta(u).
$$
Combination of upper bound in (\ref{bonferroni}) and lower bound in (\ref{lowerfields}) leads to
\BQN\label{posb}
\liminf_{n\rw\IF}\frac{\mathcal{B}_{\alpha_1,\alpha_2}(x, [0,n]^2)}{n^2}\geq \frac{\mathcal{B}_{\alpha_1,\alpha_2}(x, [0,n_1]^2)}{4n_1^2}-e^{-\mathbb{C}_1n_1^{\mathbb{C}}}.
\EQN
For $n_1>\sqrt{x}$
\BQNY
\mathcal{B}_{\alpha_1,\alpha_2}(x, [0,n_1]^2)&=&\int_{\mathbb{R}}\pk{\int_{[0,n_1]^2}\I\left( \sum_{i=1}^2 (\sqrt{2}B_{\alpha_i}(t_i) -\abs{t_i}^{\alpha_i})>s\right) dt>x}e^sds\\
&\geq&\int_\R  \pk{ \inf_{t\in [0,n_1]^2} \sum_{i=1}^2 \left(\sqrt{2}B_{\alpha_i}(t_i) -\abs{t_i}^{\alpha_i}\right)>s } e^{s} ds>0,
\EQNY
which combined with the monotonicity of $\mathcal{B}_{\alpha_1,\alpha_2}(x, [0,n_1]^2)$ in $n_1$ and (\ref{posb}) implies that for  sufficiently large $n_1$
$$\liminf_{n\rw\IF}\frac{\mathcal{B}_{\alpha_1,\alpha_2}(x, [0,n]^2)}{n^2}\geq\frac{\mathcal{B}_{\alpha_1,\alpha_2}(x, [0,n_1]^2)-4n_1^2 e^{-\mathbb{C}_1n_1^{\mathbb{C}}}}{4n_1^2}>0,$$
establishing the proof of (i). 
\\

{\underline{ \it Ad (ii)}.}
We follow notation introduced in the proof of Proposition \ref{Onepoint} for the case $\alpha_1=\beta_1$ and $\alpha_2<\beta_2$.
Let next for $u>0$
$$E_2(u)\coloneqq\left[-\left(\frac{e_u^{-1/4}\wedge\ln u}{u}\right)^{2/\beta_1},\left(\frac{e_u^{-1/4}\wedge\ln u}{u}\right)^{2/\beta_1}\right]\times \left[-\left(\frac{e_u^{-1/4}\wedge\ln u}{u}\right)^{2/\beta_2},\left(\frac{e_u^{-1/4}\wedge\ln u}{u}\right)^{2/\beta_2}\right], $$
$$ I_{k_1,k_2}(u,n)\coloneqq [k_1v_1(u)n,(k_1+1)v_1(u)n]\times[k_2v_2(u)n,(k_2+1)v_2(u)n], $$
$$ \Theta_1(u)\coloneqq 2\hat{\Gamma}(1/\beta_2+1)a_2^{1/\alpha_2}b_2^{-1/\beta_2}u^{2/\alpha_2-2/\beta_2}\Psi(u),$$
where $\hat{\Gamma}(\cdot)$ is the gamma function and
$$e_u=\sup_{0<|t_i|<\left(\frac{\ln u}{u}\right)^{2/\beta_i}, i=1,2} |e(t)|,\ \ ~e(t)=\frac{1-{\sigma}(t)}{\sum_{i=1}^2b_i|t_i|^{\beta_i}}-1, |t|\neq 0.$$
Observe that
\BQNY
\pk{\int_{E_{2}(u)}\mathbb{I}(X(t)>u)dt >\vv(u)x}&\geq& \pk{\int_{E_{1}(u,n)}\mathbb{I}(X(t)>u)dt >\vv(u)x},\\
\pk{\int_{E_{2}(u)}\mathbb{I}(X(t)>u)dt >\vv(u)x}&\leq& \pk{\int_{\bigcup_{|k_2|\leq N_2'(u,n)+1}\hat{I}_{k_2}(u,n)}\mathbb{I}(X(t)>u)dt >\vv(u)x}\\
&&\quad +\pk{\sup_{E_2(u)\setminus (\bigcup_{|k_2|\leq N_2'(u,n)+1}\hat{I}_{k_2}(u,n))}X(t)>u}.
\EQNY
Hence it follows that
\BQN\label{bonferroni1}
\Sigma_3^-(u, n_1)-\Sigma\Sigma_3(u,n_1)\leq \pk{\int_{E_{2}(u)}\mathbb{I}(X(t)>u)dt >\vv(u)x}\leq \Sigma_3^+(u,n)+\Sigma_3'(u,n),
\EQN
with
\BQNY
\Sigma_3^{\pm}(u,n)=\sum_{|k_2|\leq N_2'(u,n)\pm 1}\pk{\int_{\hat{I}_{k_2}(u,n)}\mathbb{I}(X(t)>u)dt >\vv(u)x},
\EQNY
where $I_{k_1,k_2}(u,n)$ is defined in (\ref{Ik}) and $\Sigma_3'$ and $\Sigma\Sigma_3$ are given in (\ref{eq3}) and (\ref{eq4}) respectively.
Noting that (\ref{uniform2}) also holds for $|k_2|\leq N_2'(u,n)+ 1$, we have for $x\geq 0$
\BQNY\Sigma_3^\pm(u,n)&\sim& \mathcal{B}_{\alpha_1,\alpha_2}^{a_1^{-1}b_1|t_1|^{\alpha_1},0}(x, [-n,n]\times[0,n])\sum_{|k_2|\leq N_2'(u,n)+1}\Psi(u_{k_2,n}^{\pm})\\
&\sim& \mathcal{B}_{\alpha_1,\alpha_2}^{a_1^{-1}b_1|t_1|^{\alpha_1},0}(x, [-n,n]\times[0,n])\Psi(u)\sum_{|k_2|\leq N_2'(u,n)+1}e^{-u^2b_2(|k_2|v_2(u)n)^{\beta_2}}\\
&\sim&\frac{\mathcal{B}_{\alpha_1,\alpha_2}^{a_1^{-1}b_1|t_1|^{\alpha_1},0}(x, [-n,n]\times[0,n])}{n}\Theta_1(u), \quad u\rw\IF.\EQNY
In light of Lemma \ref{gaussianfields1}, we have that for $u$ and $n$ sufficiently large
$$\Sigma\Sigma_3(u,n)+\Sigma_3'(u,n)\leq \left(\frac{\CC_2}{\sqrt{n}}+e^{-\mathbb{C}_1n^{\mathbb{C}}}\right)\Theta_1(u).$$
Dividing both sides of (\ref{bonferroni1}) by $\Theta_1(u)$ respectively and letting $u\rw\IF$, we have that
$$\frac{\mathcal{B}_{\alpha_1,\alpha_2}^{a_1^{-1}b_1|t_1|^{\alpha_1},0}(x, [-n_1,n_1]\times[0,n_1])}{n_1}-\frac{\CC_2}{\sqrt{n_1}}-e^{-\mathbb{C}_1n_1^{\mathbb{C}}} \leq \frac{\mathcal{B}_{\alpha_1,\alpha_2}^{a_1^{-1}b_1|t_1|^{\alpha_1},0}(x, [-n,n]\times[0,n])}{n}+\frac{\CC_2}{\sqrt{n}}+e^{-\mathbb{C}_1n^{\mathbb{C}}},$$
which gives that
\BQNY
\liminf_{n\rw\IF}\frac{\mathcal{B}_{\alpha_1,\alpha_2}^{a_1^{-1}b_1|t_1|^{\alpha_1},0}(x, [-n,n]\times[0,n])}{n}=\limsup_{n\rw\IF}\frac{\mathcal{B}_{\alpha_1,\alpha_2}^{a_1^{-1}b_1|t_1|^{\alpha_1},0}(x, [-n,n]\times[0,n])}{n}<\IF.
\EQNY
Moreover, we have
\BQNY
\pk{\int_{E_{2}(u)}\mathbb{I}(X(t)>u)dt >\vv(u)x}\geq \Sigma_4(u,n)-\Sigma\Sigma_4(u,n),
\EQNY
where $\Sigma\Sigma_4(u,n)$ is defined in (\ref{ssigma4}) and
\BQNY
\Sigma_4(u,n)=\sum_{|2k_2|\leq N_2'(u,n)- 1}\pk{\int_{\hat{I}_{2k_2}(u,n)}\mathbb{I}(X(t)>u)dt >\vv(u)x}.
\EQNY
By  (\ref{uniform2}),  for $x\geq 0$ we have
\BQNY\Sigma_4(u,n)&\sim& \mathcal{B}_{\alpha_1,\alpha_2}^{a_1^{-1}b_1|t_1|^{\alpha_1},0}(x, [-n,n]\times[0,n])\sum_{|2k_2|\leq N_2'(u,n)-1}\Psi(u_{k_2,n}^{-})\\
&\sim&\frac{\mathcal{B}_{\alpha_1,\alpha_2}^{a_1^{-1}b_1|t_1|^{\alpha_1},0}(x, [-n,n]\times[0,n])}{2n}\Theta_1(u), \quad u\rw\IF.\EQNY
By Lemma \ref{gaussianfields1}, for $u$ and $n$ sufficiently large, we have
$$\Sigma\Sigma_4(u,n)\leq e^{-\mathbb{C}_1n^{\mathbb{C}}}\Theta_1(u).$$
In view of (\ref{bonferroni1}) for the upper bound, we have
\BQNY
\liminf_{n\rw\IF}\frac{\mathcal{B}_{\alpha_1,\alpha_2}^{a_1^{-1}b_1|t_1|^{\alpha_1},0}(x, [-n,n]\times[0,n])}{n}&\geq&\frac{\mathcal{B}_{\alpha_1,\alpha_2}^{a_1^{-1}b_1|t_1|^{\alpha_1},0}(x, [-n_1,n_1]\times[0,n_1])}{n_1}-e^{-\mathbb{C}_1n_1^{\mathbb{C}}}.
\EQNY
Noting that for $n>\sqrt{x}$
\BQNY
\mathcal{B}_{\alpha_1,\alpha_2}^{a_1^{-1}b_1|t_1|^{\alpha_1},0}(x, [-n,n]\times[0,n])&=&\int_{\mathbb{R}}\pk{\int_{ [-n,n]\times[0,n]}\I\left( \sum_{i=1}^2 (B_{\alpha_i}(t_i) -\abs{t_i}^{\alpha_i})-a_1^{-1}b_1|t_1|^{\alpha_1}>s\right) dt>x}e^sds\\
&\geq&\int_\R  \pk{ \inf_{t\in [-n,n]\times[0,n]} \left(\sum_{i=1}^2 \left(B_{\alpha_i}(t_i) -\abs{t_i}^{\alpha_i}\right)-a_1^{-1}b_1|t_1|^{\alpha_1}\right)>s } e^{s} ds>0,
\EQNY and by the monotonicity of $\mathcal{B}_{\alpha_1,\alpha_2}^{a_1^{-1}b_1|t_1|^{\alpha_1},0}(x, [-n,n]\times[0,n])$     with respect to $n$, we have, for $n_1$ sufficiently large,
\BQNY\liminf_{n\rw\IF}\frac{\mathcal{B}_{\alpha_1,\alpha_2}^{a_1^{-1}b_1|t_1|^{\alpha_1},0}(x, [-n,n]\times[0,n])}{n}&\geq&\frac{\mathcal{B}_{\alpha_1,\alpha_2}^{a_1^{-1}b_1|t_1|^{\alpha_1},0}(x, [-n_1,n_1]\times[0,n_1])}{n_1}-e^{-\mathbb{C}_1n_1^{\mathbb{C}}}>0.
\EQNY
This completes the proof of (ii).
\\

{\underline{ \it Ad (iii)}.}
We follow notation introduced in the proof of Proposition \ref{Onepoint} for the case $\alpha_i=\beta_i$, $i=1,2$
Observe that
 \BQN\label{sigma3}
\Sigma_5(u,n)\leq \pk{\int_{E'(u,n)}\mathbb{I}(X(t)>u)dt >\vv(u)x}\leq \Sigma_5(u,n)+\Sigma\Sigma_5(u,n),
 \EQN
 where $E'(u,n)=\bigcup_{(k_1,k_2)\in K_{u,n}}I_{k_1,k_2}(u,n)$ and
 $$\Sigma_5(u,n)=\pk{\int_{\hat{I}(u,n)}\mathbb{I}(X(t)>u)dt >\vv(u)x}, $$ $$
  \Sigma\Sigma_5(u,n)=\sum_{|k_i|\leq N_i'(u,n), k_i\neq -1,0, i=1,2}\pk{\sup_{t\in I_{k_1,k_2}(u,n)}\overline{X}(t)>u_{n,k_1,k_2}^{-}},$$
 with $u_{n,k_1,k_2}^{-}$ defined in (\ref{uk}) and $\hat{I}(u,n)$ in (\ref{EE1}). In light of (\ref{uniform3}) and (\ref{var}), we have that for $u$ sufficiently large
 \BQNY
 \Sigma\Sigma_5(u,n)&\leq& \mathcal{B}_{\alpha_1,\alpha_2}(x, [0,n]^2)\sum_{|k_i|\leq N_i'(u,n), k_i\neq -1,0, i=1,2}\Psi(u_{n,k_1,k_2}^{-})\\
 &\leq&  \mathcal{B}_{\alpha_1,\alpha_2}(x, [0,n]^2)\Psi(u)\sum_{|k_i|\leq N_i'(u,n), k_i\neq -1,0, i=1,2}e^{-a_1^{-1}b_1|k_1^*n|^{\beta_1}-a_2^{-1}b_2|k_2^*n|^{\beta_2}}\\
 &\leq& \mathcal{B}_{\alpha_1,\alpha_2}(x, [0,n]^2)e^{-Q_1(n^{\beta_1}+n^{\beta_2})}\Psi(u),
 \EQNY
 where
 $k_i^*=k_iI_{\{k_i>0\}}+(|k_i|-1)I_{\{k_i<0\}}, i=1,2.$

 Hence dividing (\ref{sigma3}) by $\Psi(u)$  and letting $u\rw\IF$, we have for any $n,n_1>\sqrt{x}$
 \BQNY
 0<\mathcal{B}_{\alpha_1,\alpha_2}^{a_1^{-1}b_1|t_1|^{\alpha_1},a_2^{-1}b_2|t_2|^{\alpha_2}}(x, [-n,n]^2)\leq \mathcal{B}_{\alpha_1,\alpha_2}^{a_1^{-1}b_1|t_1|^{\alpha_1},a_2^{-1}b_2|t_2|^{\alpha_2}}(x, [-n_1,n_1]^2)+\mathcal{B}_{\alpha_1,\alpha_2}(x, [0,n_1]^2)e^{-Q_1(n_1^{\beta_1}+n_1^{\beta_2})}.
 \EQNY
 Letting $n\rw\IF$ with $n_1$ fixed in the above inequality, we complete the proof.
 \QED

{\bf Proof of \eqref{upper-lower}}: Observe that
\BQNY
\frac{\Psi(u_{n,k_1,k_2}^{-})}{\Psi(u_{n,k_1,k_2}^{+})}\sim e^{\frac{\left(u_{n,k_1,k_2}^{+}\right)^2-\left(u_{n,k_1,k_2}^{-}\right)^2}{2}}, \quad u\rw\IF
\EQNY
uniformly with respect to $0\leq |k_i|\leq N_i{'}(u,n), i=1,2$. Furthermore, by (\ref{var}), for $u$ sufficiently large
\BQNY
\left(u_{n,k_1,k_2}^{+}\right)^2-\left(u_{n,k_1,k_2}^{-}\right)^2&=&u^2\left(\sup_{t\in I_{k_1,k_2}(u,n)}\frac{1}
{{\sigma^2(t)}}-\inf_{t\in I_{k_1,k_2}(u,n)}\frac{1}{{\sigma^2(t)}}\right)\\
&=&u^2\sup_{s, t\in I_{k_1,k_2}(u,n)}\frac{{\sigma^2(t)}-{\sigma^2(s)}}{{\sigma^2(t)}{\sigma^2(s)}}\\
&\leq &4u^2\sup_{s, t\in I_{k_1,k_2}(u,n)}|{\sigma(t)}-{\sigma(s)}|\\
&=& 4u^2 \sup_{s, t\in I_{k_1,k_2}(u,n)}\left|(1+e(t))\sum_{i=1}^2b_i|t_i|^{\beta_i}-(1+e(s))\sum_{i=1}^2b_i|s_i|^{\beta_i}\right|\\
&\leq& 4u^2 \sup_{s, t\in I_{k_1,k_2}(u,n)}\left|\sum_{i=1}^2b_i|t_i|^{\beta_i}-\sum_{i=1}^2b_i|s_i|^{\beta_i}\right|+8u^2 \sup_{ t\in I_{k_1,k_2}(u,n)}|e(t)|\sum_{i=1}^2b_i|t_i|^{\beta_i}\\
&\leq& 4u^2 \sum_{i=1}^2b_i\beta_i|\theta_i|^{\beta_i-1}v_i(u)n+8u^2 \sup_{ t\in I_{k_1,k_2}(u,n)}|e(t)|\sum_{i=1}^2b_i|t_i|^{\beta_i},
\EQNY
where $e(t)=\frac{1-{\sigma(t)}}{\sum_{i=1}^2b_i|t_i|^{\beta_i}}-1, |t|\neq 0$ and $\theta_i\in (k_iv_i(u)n, (k_i+1)v_i(u)n)$. Using the fact that
$$N_i'(u,n)=\left[\frac{(e_u^{-1/4}\wedge\ln u)^{2/\beta_i}}{u^{2/\beta_i}v_i(u)n}\right]~\text{and}~\lim_{u\rw\IF}e_u=0,$$
we have that
$$u^2 \sup_{ t\in I_{k_1,k_2}(u,n)}|e(t)|\sum_{i=1}^2b_i|t_i|^{\beta_i}\leq 2e_u\sum_{i=1}^2b_i(e_u^{-1/4}\wedge\ln u)^2\rw 0,$$
as $u\rw\IF$ uniformly with respect to $0\leq |k_i|\leq N_i{'}(u,n), i=1,2$. For $\beta_i\geq 1, i=1,2,$ 
\BQNY
u^2 
\sum_{i=1}^2b_i\beta_i|\theta_i|^{\beta_i-1}v_i(u)n
&\leq& u^2 \sum_{i=1}^2b_i\beta_i\left(\frac{\ln u}{u}\right)^{\frac{2(\beta_i-1)}{\beta_i}}v_i(u)n\\
&\leq&\sum_{i=1}^2 2a_i^{-1/\alpha_i}b_i\beta_i u^{2/\beta_i-2/\alpha_i}(\ln u)^{\frac{2(\beta_i-1)}{\beta_i}}n\rw 0, \quad u\rw\IF
\EQNY
uniformly with respect to $0\leq |k_i|\leq N_i{'}(u,n), i=1,2$, where $(\theta_1,\theta_2)\in I_{k_1,k_2}(u,n)$.
For $0<\beta_i<1, i=1,2,$ 
\BQNY
u^2 \sup_{s, t\in I_{k_1,k_2}(u,n)}\left|\sum_{i=1}^2b_i|t_i|^{\beta_i}-\sum_{i=1}^2b_i|s_i|^{\beta_i}\right|
&\leq& u^2 
\sum_{i=1}^2b_i\beta_i|\theta_i|^{\beta_i-1}v_i(u)n\\
&\leq& u^2 \sum_{i=1}^2b_i\beta_i|v_i(u)n|^{\beta_i}\rw 0, \quad u\rw\IF,
\EQNY
holds uniformly for $0\leq |k_i|\leq N_i{'}(u,n), k_i\neq -1,0,  i=1,2.$ For $0<\beta_i<1,  k_i=-1,0,i=1,2$
\BQNY
u^2 \sup_{s, t\in I_{k_1,k_2}(u,n)}\left|\sum_{i=1}^2b_i|t_i|^{\beta_i}-\sum_{i=1}^2b_i|s_i|^{\beta_i}\right|
&\leq& u^2 \sup_{s, t\in I_{k_1,k_2}(u,n)}\left(\sum_{i=1}^2b_i|t_i|^{\beta_i}+\sum_{i=1}^2b_i|s_i|^{\beta_i}\right)\\
&\leq& 2u^2\sum_{i=1}^2b_i|v_i(u)n|^{\beta_i}\\
&=& 2\sum_{i=1}^2a_i^{-\beta_i/\alpha_i}b_in^{\beta_i}u^{2-2\beta_i/\alpha_i}\rw 0, \quad u\rw\IF.
\EQNY
Therefore, we can conclude that
\BQNY
\left(u_{n,k_1,k_2}^{+}\right)^2-\left(u_{n,k_1,k_2}^{-}\right)^2\rw 0
\EQNY
as $u\rw\IF$ uniformly with respect to $0\leq |k_i|\leq N_i'(u,n), i=1,2$ establishing the proof. \QED


{\bf Acknowledgement}: K. D\c ebicki
was partially supported by NCN Grant No 2018/31/B/ST1/00370 (2019-2022).
Partial financial support  from  Swiss National Science Foundation Grant 200021-196888 is kindly acknowledged.  Additionally, this project is financed by the Ministry of Science and Higher Education in Poland under the programme "Regional Initiative of Excellence" 2019 - 2022 project number 015/RID/2018/19.

\bibliographystyle{plain}

 \bibliography{queue2cd}

\end{document}